\documentclass[11pt]{article}
\usepackage{epsfig}
\usepackage{verbatim}
\usepackage{amsmath}
\usepackage[pdfmenubar=false, letterpaper=false, pdfpagelabels=true]{hyperref}
\usepackage{makeidx,showidx}
\usepackage{pst-node}  

\begin{document}

\newcommand{\bm}[1]{\mbox{\boldmath$#1$}}
\def\mvec#1{{\bm{#1}}}   % vector in math mode; bold

\title{A defense of Columbo\\
(and of the use of Bayesian inference in forensics)\\
-- A multilevel introduction to probabilistic reasoning --
}
\author{G.~D'Agostini \\
Universit\`a ``La Sapienza'' and INFN, Roma, Italia \\
{\small (giulio.dagostini@roma1.infn.it,
 \url{http://www.roma1.infn.it/~dagos})}
}

\date{}

\maketitle

\begin{abstract}
Triggered by a recent interesting New Scientist article
on the too frequent incorrect use of probabilistic 
evidence in courts, I introduce the basic concepts of 
probabilistic inference with a toy model, and
discuss several important issues that need to be
understood in order
to extend the basic reasoning to real life cases.
In particular, I emphasize
the often neglected point
that  degrees of beliefs
are updated not by `bare facts' alone, but by all
available information pertaining to them, 
including how they have been acquired. 
In this light I show that,
contrary to what claimed in that article, 
there was no ``probabilistic pitfall''
in the Columbo's episode pointed as example
of ``bad mathematics'' yielding ``rough justice''. 
Instead, such a criticism could have a `negative
reaction' to the article itself and to the
use of Bayesian reasoning in courts,
 as well as in all other places
in which probabilities need to be assessed and 
decisions need to be made. Anyway, besides 
introductory/recreational aspects, the paper touches important 
questions, like: role and evaluation of priors;
subjective evaluation of Bayes factors; 
role and limits of intuition; 
`weights of evidence' and `intensities of beliefs'
(following Peirce) and `judgments leaning' 
(here introduced), including their uncertainties
and combinations;
role of relative frequencies to assess and express
beliefs; pitfalls due to `standard' statistical education;
weight of evidences mediated by  testimonies. 
A small introduction to Bayesian networks, based 
on the same toy model (complicated by the 
possibility of incorrect testimonies) and implemented 
using Hugin software, is also 
provided, to stress the importance of formal, computer aided
probabilistic reasoning.
\end{abstract}
\vspace{1.2cm}
{\small
\begin{flushright}
{\sl ``Use enough common sense to know } \\
{\sl when ordinary common sense does not apply'' } \\
{(I.J. Good's {\it guiding principle of all science})}\\
\end{flushright}
}

\newpage

\section{Introduction}
A recent New Scientist article~\cite{NS}
deals with errors in courts due to ``bad mathematics'',
advocating the use of the so-called Bayesian methods to avoid them.
Although most examples of resulting ``rough justice'' come
from real life cases, the first
``probabilistic pitfall'' 
is taken from crime fiction, 
namely from a ``1974 episode of the cult
US television series'' {\it Columbo}, in which 
a ``society photographer has killed his wife and disguised 
it as a bungled kidnapping.'' 

The pretended mistake happens in the 
concluding scene,
when ``the hangdog detective 
[\ldots]
induces the murderer to grab from a shelf of 12 cameras 
the exact one used to snap the victim before she was killed.''
According to the article author (or to experts on which 
scientific journalists often rely on) the question is that 
``killer or not, anyone
would have a 1 in 12 chance of picking the same 
camera at random. That kind of evidence would 
never stand up in court.'' Then a sad doubt is raised,
``Or would it? In fact, such probabilistic pitfalls 
are not limited to crime fiction.''

\begin{figure}
\centering\epsfig{file=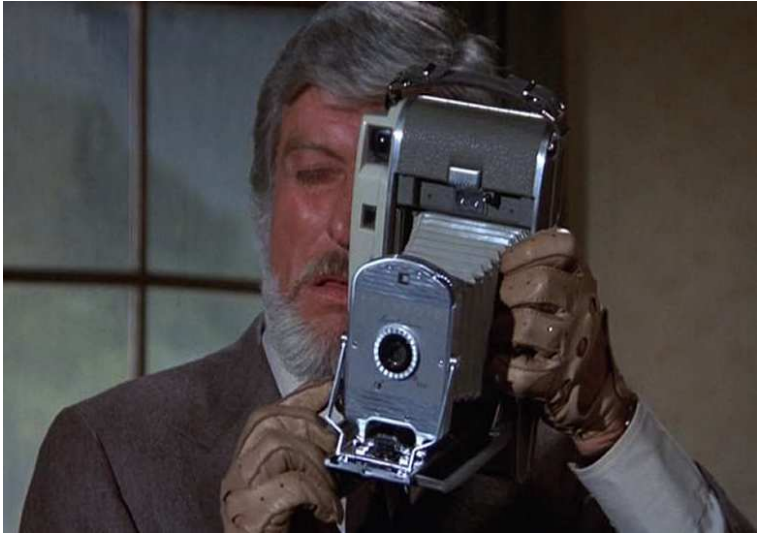,clip=}
\caption{\small \sf Peter Galesco, played by Dick Van Dyke,
taking the picture of its wife before killing her.}
\label{fig:camera}
\end{figure}

Being myself not particularly fond of this kind of entertainment
(perhaps with a little exception of the Columbo series,
that I watch casually),
I cannot tell
how much crime fiction literature and movies are affected by 
``probabilistic
pitfalls''. Instead, I can give firm witness
that scientific practice 
is plenty of mistakes of the kind reported 
in Ref.\cite{NS}, that happen even
in fields the general public 
would hardly suspect, like frontier physics, 
whose protagonists are supposed to have a skill
in mathematics superior to police officers and lawyers.

But it is not just a question of math skill (complex 
calculations are usually done without mistakes), 
but of {\it probabilistic reasoning} (what to calculate!). 
This is a quite old story.
In fact, as David Hume complained 260 years ago~\cite{DH},
\begin{quote}
{\sl 
"The celebrated Monsieur Leibniz has observed
it to be a defect in the
common systems of logic, that they are very copious 
when they explain the
operations of the understanding in the forming of demonstrations,
 but are
too concise when they treat of probabilities, and 
those other measures of
evidence on which life and action entirely depend, 
and which are our guides
even in most of our philosophical speculations."
}
\end{quote}
It seems to me that the general situation has not improved much. 
Yes, `statistics' (a name that, meaning too much, risks 
to mean little) is taught in colleges and universities
to students of several fields, but distorted
by the `frequentistic approach', according to which 
one is not allowed to speak
of probabilities of causes.
This is, in my opinion,
the original sin that gives grounds for
a large number of probabilistic mistakes
even by otherwise very valuable scientists and practitioners
(see e.g. 
chapter 1 of Ref.~\cite{BR}).

Going back to the ``shambling sleuth Columbo'',
being my wife and my daughter 
his fans, it happens
we own the DVD collections of the first seven seasons.
It occurred then I watched 
with them, 
not much time ago (perhaps last winter),
the `incriminated', superb 
episode {\it Negative Reaction}~\cite{NR}, one of the best
performances of Peter Falk playing the role
of the famous lieutenant. 
However, reading the mentioned New Scientist article,
I did not remember I had a `negative
reaction' from the final scene, although 
I use and teach Bayesian methods for a large
variety of applications. Did I overlook something?

I watched again the episode and I was again convinced
Columbo's last move was a conclusive 
checkmate.\footnote{Just writing this note, 
I have realized that
the final scene is directed so well
that, not only the way the photographer 
loses control and commits his fatal mistake looks very credible, 
but also spectators forget he could play valid countermoves, 
not depending
on the negative of the pretended destroyed picture 
(see footnote \ref{fn:nota}). 
Therefore, rather than chess, the name of 
the game is {\it poker}, and 
Columbo's bluff is able to induce
the murderer to provide a crucial piece of evidence
to finally incriminate him.}
Then I have 
invited some friends, all with physics or mathematics 
degree and somewhat knowledgeable of the Bayesian 
approach,
to enjoy an evening together during the recent end of year 
holidays in order to let them make up their minds
whether Columbo
had good reasons to take Paul Galesco, magnificently
impersonated by Dick Van Dyke, in front of the court
(Bayes or not, we had some fun\ldots).

The verdict was unanimous: Columbo was fully absolved
or, more precisely, there was nothing to reproach 
the story writer, Peter S. Fischer. The convivial 
after dinner jury also requested me to write
a note on the question, possibly with a short, 
self-contained introduction to the `required math'.   
Not only to `defend Columbo' or, more properly,
his writer, but, and more seriously, 
to defend the Bayesian approach, and in particular
its applications in forensic science. 
In fact, we all deemed the beginning 
paragraphs of the New Scientist article could throw a bad
light on the rest of the contents. 

Imagine a casual reader of the article, possibly 
a lawyer, a judge or a student in forensic science, 
to which the article was virtually
addressed, and who might have seen {\it Negative Reaction}.
Most likely he/she considered legitimate 
the charges of the policemen 
against the photographer. 
The `negative reaction'
would be that the reader would consider the rest
of the article a support of dubious validity 
to some `strange math'
that can never substitute the human intuition
in a trial.\footnote{This kind of objection,
in defense of what is often nothing but 
``the capricious ipse dixit of authority''\cite{JHN},
from which we should instead ``emancipate''\cite{JHN},
is quite frequent.
It is raised not only by judges, who tend to claim their job is
"to evaluate evidence not by means of a formula... 
but by the joint application of their 
individual common sense."\cite{NS}, 
but also by other categories of people 
who take important decisions, like
doctors, managers and politicians.
\label{note:common-sense}} 
Not a good service to the `Bayesian cause'. 
(Imagine somebody trying to convince you with arguments
you hardly understand and who begins
asserting something you consider manifestly false.) 

In the following section I introduce the 
basic elements of Bayesian reasoning
(subsection \ref{ss:JL} can be skipped on first reading), 
using a toy model as guiding example  
in which the  analysis of ref.~\cite{NS} (``1 in 12'', or, more
precisely ``1 in 13'') holds. 
Section \ref{sec:Columbo_priors}
shows how such a kind of evidence would change
Columbo's and jury's opinion.  
Then I discuss in section \ref{sec:finale}
why a similar argument
does not apply to the clip in which 
Columbo finally frames Galesco, and 
why all witnesses of the crucial actions 
(including TV watchers, with the exception
of the author of Ref.~\cite{NS} and perhaps a few others)  
and an hypothetical 
court jury
(provided the scene had been properly reported)
had to be absolutely positive the photographer killed
his wife (or at least 
he knew who did it in his place).

The rest of the paper might be marginal, if you
are just curious to know why I have a different opinion than
Ref.~\cite{NS}, although I agree
on the validity of Bayesian reasoning. 
In fact, at the end of the work,
this paper is not the `short note' 
initially planned. The reason is that the past months
I had many discussions on some of the questions
treated here with people from several fields. 
I have realized once more that it is not easy
to put the basic principles at work 
if some important issues are not well understood.
People are used to solving their
statistical problems with `ad hoc' formulae (see Appendix H) and therefore
tend to add some `Bayesian recipes' in their 
formularium. 
It is then too high the risk that one looks
at simplified methods -- Bayesian methods require a bit more
thinking and computation that others! -- that are even
advertised as `objective'. 
Or one just refuses to use any math, on the defense
of pure intuition. (By the way, this is an important point and
I will take the opportunity to comment 
on the apparent contradictions between 
intuition and formal evaluation of beliefs,
defending \ldots both,
but encouraging the
use of the latter, superior to the former 
in complex situations -- see in particular Appendix C).

So, to conclude the introduction, this document
offers several levels
of reading:
\begin{itemize}
\item If you are only interested to Columbo's story, 
      you can just jump straight to section \ref{sec:finale}.
\item If you also (or regardless of Columbo) want to have an opportunity
      to learn the basic rules of Bayesian inference,
      subsections \ref{ss:Bayes_theorem},  \ref{ss:Bayes_theorem_2}
      and \ref{ss:many_pieces}, based on a simple master 
      example, have been written on the purpose. 
      Then you might appreciate the advantage of logarithmic
      updating (section \ref{ss:JL}) and perhaps see how it applies to 
      the AIDS example of Appendix F.
      \\ \mbox{}
      
\item If you already know the basics of the probabilistic
      reasoning, but you wonder how it can be applied
      into real cases, 
      then section \ref{sec:real_life} should help, together
      with some of the appendices.
\item If none of the previous cases is yours 
      (you might
      even be an expert of the field), 
      you can simply browse the document. Perhaps
      some appendices  or 
      subsections might 
      still be of your interest.
\item Finally, there is the question of the many footnotes,
      which can break the pace of the reading. 
      They are not meant to  be necessarily 
      read sequentially along with the main text and could be skipped
      on a first fast reading (in fact, this document 
      is closer to an hypertext than to a standard article.)
\end{itemize}

Enjoy!
\\ \break \mbox{}

\section{One in thirteen -- Bayesian reasoning illustrated
with a toy model}\label{sec:1in13}
Let us  leave aside Columbo's cameras  for a while
and begin with a different, simpler, 
stereotyped situation easier to analyze.

Imagine there are two types of boxes, $B_1$, that only contain
white balls ($W$), and $B_2$, that contain one
white balls and twelve black 
(incidentally, just to be precise,
although the detail is absolutely irrelevant, 
we have to infer from Columbo's words,
{\it ``You didn't touch any
of these twelve cameras. You picked up that one''},
the cameras  were thirteen).

You take at random a box and extract a ball.
The resulting color is
{\it white}.
You might be interested to evaluate 
the probability that the box is of type $B_1$,
in the sense of stating in a quantitative way
how much you believe
this hypothesis.
In formal terms we are interested in $P(B_1\,|\,W,I)$, 
knowing that $P(W\,|\,B_1,I)=1$ and 
$P(W\,|\,B_2,I)=1/13$, a problem that can be sketched as
\begin{eqnarray}
\left\{\!\!
\begin{array}{rcl} P(W\,|\,B_1,I)& = &1 \\
P(W\,|\,B_2,I)&=&{1}/{13}
\end{array}
\right.
 \hspace{0.8cm}
&\Rightarrow&\hspace{0.3cm}
 P(B_1\,|\,W,I)=\, \bm{?}
\end{eqnarray}
[Here `$|$' stands for `given', or `conditioned by';
$I$ is the general (`background') status of 
information under which this probability is assessed;
`$W,I$' or `$B_i,I$' after `$|$' indicates that both conditions
are relevant for the evaluation of the probability.]

A typical mistake
at this point is to
confuse $P(B_1\,|\,W,I)$ with $P(W\,|\,B_1,I)$,
or, more often, 
$P(B_2\,|\,W,I)$ with $P(W\,|\,B_2,I)$,
as largely discussed in Ref.~\cite{NS}.
Hence we need to learn how 
to turn properly $P(W\,|\,B_1,I)$ into $P(B_1\,|\,W,I)$ 
 using the rules of probability 
theory.

\subsection{Bayes theorem and Bayes factor}\label{ss:Bayes_theorem}
The `probabilistic inversion' 
$P(W\,|\,B_1,I) \rightarrow P(B_1\,|\,W,I)$
can {\it only}\footnote{Beware of methods 
that provide `levels of confidence', or something 
like that, without using Bayes' theorem! See also 
footnote \ref{fn:nopriors} and Appendix H.} be performed
using
the so-called {\it Bayes' theorem}, a simple consequence
of the fact that, given the {\it effect} $E$ and some
{\it hypotheses} $H_i$ concerning its possible cause,
the joint probability of $E$ and $H_i$, 
conditioned by the {\it background
information}\footnote{The background information $I$
represents all we know about the hypotheses
and the effect considered. Writing $I$ in all 
expressions could seem a pedantry, 
but it isn't. 
For example, if we would just write 
$P(E)$ in these formulae, instead of  $P(E\,|\,I)$,
 one might be tempted to take 
this probability equal to one, 
``because the observed event is a well established fact',
that has happened and is then certain. 
But it is not this certainty 
that enters these formulae, but rather the probability `that
fact could happen' 
in the light of `everything we knew'  about it (`$I$').} 
$I$, can be written as 
\begin{eqnarray}
P(E\cap H_i\,|\,I) &=& P(E\,|\,H_i,I)\cdot P(H_i\,|\,I) = 
                    P(H_i\,|\,E,I)\cdot P(E\,|\,I)\,,\label{eq:joint_prob}
\end{eqnarray}
where `$\cap$' stands for a logical `AND'. 
From the second equality of the last
equation we get%
\begin{eqnarray}
P(H_i\,|\,E,I) &=& \frac{P(E\,|\,H_i,I)}{P(E\,|\,I)}\cdot P(H_i\,|\,I)\,,
\label{eq:PHi|E}
\end{eqnarray}
that is one of the ways to express Bayes' 
theorem.\footnote{Bayes' theorem can be often found in the form
\begin{eqnarray*}
P(H_i\,|\,E,I) &=& \frac{P(E\,|\,H_i,I)\cdot  P(H_i\,|\,I)}
                        {\sum_i P(E\,|\,H_i,I)\cdot  P(H_i\,|\,I)}\,,
\end{eqnarray*}
valid if we deal with a class of incompatible hypotheses
[i.e. $P(H_i\cap H_j\,|\,I)=0$ and $\sum_i P(H_i\,|\,I)=1$]. 
In fact, in this case 
a general rule of probability theory [Eq.~(\ref{eq:rul6}) in Appendix A]
allows us
to rewrite the denominator of Eq.~(\ref{eq:PHi|E}) 
as  $\sum_i P(E\,|\,H_i,I)\cdot  P(H_i\,|\,I)$. 
In this note, dealing only with two hypotheses, 
we prefer to reason in terms of probability ratios,
as shown in Eq.~(\ref{eq:Bayes_factor}).
}%- fine footnote

Since a similar expression holds for 
any other hypothesis $H_j$, 
dividing member by member the two expressions
we can restate the theorem in terms of the 
relative beliefs, that is
\begin{eqnarray}
\underbrace{\frac{P(H_i\,|\,E,I)}{P(H_j\,|\,E,I)}}_{\mbox{updated odds}} 
%`posterior odds'
&=& 
\underbrace{\frac{P(E\,|\,H_i,I)}{P(E\,|\,H_j,I)}
           }_{\begin{array}{c} \mbox{updating factor} \\ \mbox{(`{\it Bayes factor}')} \end{array}}
 \times
\underbrace{\frac{P(H_i\,|\,I)}{P(H_j\,|\,I)}}_{\mbox{initial odds}}\,:
\label{eq:Bayes_factor}
\end{eqnarray}
the initial ratio of beliefs (`odds') is updated by the 
so-called {\it Bayes factor}, that depends on how likely
{\it each} hypothesis can produce that effect.\footnote{Note that, 
while in the case of only two hypotheses entering the inferential
game their initial probabilities are related by
$P(H_2\,|\,I) = 1- P(H_1\,|\,I)$, the probabilities of the effects
$P(E\,|\,H_1,I)$ and $P(E\,|\,H_2,I)$ have usually nothing to do 
with each other.}
Introducing $O_{i,j}$ and $\mbox{BF}_{i,j}$, 
with obvious meanings, we can rewrite
Eq.~(\ref{eq:Bayes_factor}) as 
\begin{eqnarray}
O_{i,j}(E,I) &=& \mbox{BF}_{i,j}(E,I) \times O_{i,j}(I)\,.
\label{eq:Bayes_factor_1}
\end{eqnarray}
Note that, if the initial odds are unitary, than 
the final odds are equal to the updating factor. Then, {\it Bayes factors
can be interpreted as odds due  only  to an individual
piece of evidence, \underline{if} the two hypotheses were
considered initially equally
likely}.\footnote{Those who want to base the inference only
on  the probabilities of the observations given the hypotheses,
in order to ``let the data speak themselves'', might be 
in good faith, but their noble intention
does dot save them from dire
mistakes~\cite{BR}. (See also footnotes \ref{fn:nopriors} and 
\ref{fn:Jaynes}, as well as Appendix H.)} 
This allows us to rewrite $\mbox{BF}_{i,j}(E,I)$ as 
$\tilde{O}_{i,j}(E,I)$, where the tilde is to remind that {\it they
are not properly odds}, but rather  
`{\it pseudo-odds}'. We get then an expression in which all terms
have {\it virtually} uniform meaning:
\begin{eqnarray}
O_{i,j}(E,I) &=& \tilde{O}_{i,j}(E,I) \times O_{i,j}(I)\,.
\label{eq:Bayes_factor_2}
\end{eqnarray}
If we have only two hypotheses, we get simply
%\begin{eqnarray}
$O_{1,2}(E,I) = \tilde{O}_{1,2}(E,I) \times O_{1,2}(I)$\,.
%\end{eqnarray}
If the updating factor is unitary, then the piece of
evidence does not modify our opinion on the two hypotheses
(no matter how small can numerator and denominator be,
as long as their ratio remains finite and unitary! -- see Appendix G
for an example worked out in details);
when $\tilde{O}_{1,2}(E,I)$ vanishes, then 
hypothesis $B_1$ becomes impossible 
(``{\it it is falsified}''); if  instead it is  infinite 
(i.e. the denominator vanishes), then it is the other 
hypothesis to be impossible. (The undefined case $0/0$
means that we have to look for other hypotheses 
to explain the effect.\footnote{Pieces of evidence 
modify, in general, 
relative beliefs. When we turn relative beliefs into
absolute ones in a scale ranging from 0 to 1, 
we are {\it always} making the implicit assumption
that the possible hypotheses are only those
of the class considered. 
If other hypotheses are added, the relative beliefs do
not change, while the absolute ones do. This is the reason
why an hypothesis can eventually be falsified, if 
$P(E\,|\,H_i,I)=0$, but an absolute truth, i.e. 
$P(E\,|\,H_j,I)=1$, depends on which class of hypotheses 
is considered. Stated in other words, 
in the realm of probabilistic inference 
{\it falsities can be absolute, but 
truths are always relative}.})

\subsection{Role of priors}\label{ss:Bayes_theorem_2}
Applying the updating reasoning to our box game,
the Bayes factor of interest is
\begin{eqnarray}
\tilde O_{1,2}(W,I) &=&\frac{P(W\,|\,B_1,I)}{P(W\,|\,B_2,I)}
= \frac{1}{1/13}= 13\,.
\end{eqnarray}
As it was remarked, this number would give the required odds
{\it if} the hypotheses were initially equally likely. 
But how strong are the initial relative beliefs 
on the two hypotheses? `Unfortunately', 
we cannot perform a probabilistic inversion
if we are unable to assign somehow {\it prior probabilities} to the
hypotheses we are interested in.\footnote{You might
be reluctant to adopt this way of reasoning, 
objecting
``I am unable to state priors!'', or
``I don't want to be influenced by prior!'', or even
``I don't want to state degrees of beliefs, but only 
real probabilities''. No problem, provided you stay
away from probabilistic inference (for example you 
can enjoy fishing or hiking -- but I hope you are aware of the
large amount of prior beliefs 
involved in these activities too!).
Here I can only advice you, provided you are interested in 
evaluating probabilities
of `causes' from effects, not to overlook prior information
and not to  blindly trust statistical methods and 
software packages advertised as prior-free, 
unless you don't want to risk to arrive at very bad
conclusions. 
For more comments on the question see Ref.~\cite{BR}, 
footnote \ref{fn:Jaynes} and Appendix H.
\label{fn:nopriors}
} 
\begin{figure}
\centering\epsfig{file=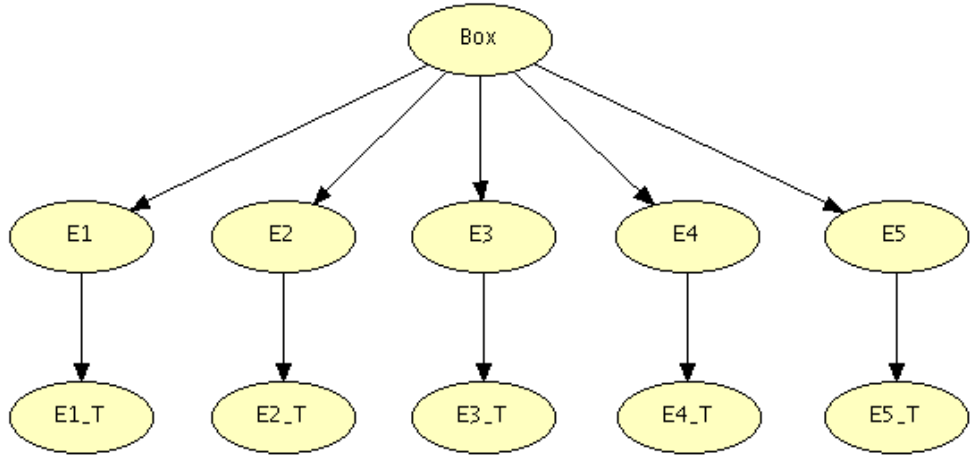,clip=,}
\caption{\small \sf Graphical representation of the causal connections
Box$\rightarrow E_i$, where $E_i$ are the effects 
(White/Black at each extraction). These effects are then causes
of other effects ($E_i\_T$), the reported colors, where `$T$'
stands for  `testimony'.  
The arrows connecting the various `nodes' represent
conditional probabilities. The model will be fully exploited
in Appendix J.}
\label{fig:bn}
\end{figure}
Indeed, in the formulation of the problem I  on purpose
passed over the relevant pieces of information
to evaluate the prior probabilities (it was said that
``there are two types of boxes'', 
not ``there are two boxes''!). If we specify that 
we had $n_1$ boxes of type $B_1$
and $n_2$ of the other kind, 
then the initial odds
are $n_1/n_2$ and the final ones will be 
\begin{eqnarray}
O_{1,2}(W,I) &=& \tilde{O}_{1,2}(W,I)\times O_{1,2}(I) \\
            &=& 13\times  \frac{n_1}{n_2},
\end{eqnarray}
from which we get (just requiring that the probability of the two
hypotheses have to sum up to one\footnote{
If $H_1$ and  $H_2$ are generic, complementary hypotheses 
we get, calling $b$ the Bayes factor of $H_1$ versus $H_2$
and $x_0$ the initial odds to simplify the notation,
 the following convenient expressions to evaluate the probability
of $H_1$:
\begin{eqnarray*}
P(H_1\,|\,x_0,b) &=& \frac{b\,x_0}{1+b\,x_0} = 
\frac{b}{b+1/x_0} = \frac{x_0}{x_0+1/b}\,.
\end{eqnarray*}
})
\begin{eqnarray}
P(B_1\,|\,W,I_0) &=& \frac{13}{13+n_2/n_1}\,.
\end{eqnarray}
If the two hypotheses were {\it initially} considered  equally likely,
then the evidence $W$ makes $B_1$ 13 times more believable
than $B_2$, i.e. $P(B_1\,|\,W,I_0)=13/14$, or approximately 93\%.
On the other hand, if $B_1$ was {\it a priori} much less
credible than $B_2$, 
for example by a factor 13,
just to play with round numbers, the same evidence
made $B_1$ and $B_2$ equally likely. Instead, if we were
initially in strong favor of $B_1$, considering it for instance
13 times more plausible than $B_2$, that evidence 
turned this factor into 169, making us 99.4\% confident 
-- {\it highly confident}, some would even 
say `practically sure'! --
that the box is of type $B_1$. 

\subsection{Adding pieces of evidence}\label{ss:many_pieces}
Imagine now the following variant of the previous toy
experiment. After the white ball is observed, you put it again
in the box, shake well and make a second extraction. 
You get white the second time too. Calling $W_1$ and $W_2$ the 
two observations, 
we have now:\footnote{Note that we
are still using Eq.~(\ref{eq:Bayes_factor}), although we are dealing
now with more complex events and complex hypotheses, 
logical AND of simpler ones. Moreover, Eq.~(\ref{eq:seq:2}) 
is obtained from Eq.~(\ref{eq:seq:1}) making use of 
the formula (\ref{eq:joint_prob}) of joint probability,
that gives $P(W_1,W_2\,|\,B_1,I) = P(W_2\,|\,W_1,B_1,I)\times
P(W_1\,|\,B_1,I)$ and an analogous formula for $B_2$. 
Note also that, going from 
Eq.~(\ref{eq:seq:2}) to Eq.~(\ref{eq:seq:3}), 
$P(W_2\,|\,W_1,B_i,I_0)$ has been rewritten as
$P(W_2\,|\,B_i,I_0)$ to emphasize that the probability
of a second white ball,
conditioned by the box composition and the 
result of the first extraction,  
depends indeed only on the box content
and not on the previous outcome (`extraction after re-introduction').
}
\begin{eqnarray}
\frac{P(B_1\,|\,W_1,W_2,I_0)}{P(B_2\,|\,W_1,W_2,I_0)} &=&
\frac{P(W_1,W_2\,|\,B_1,I_0)}{P(W1,W_2\,|\,B_2,I_0)} \times  
\frac{P(B_1\,|\,I_0)}{P(B_2\,|\,I_0)}\label{eq:seq:1} \\
&=& \frac{P(W_2\,|\,W_1,B_1,I_0)\cdot P(W_1\,|\,B_1,I_0)}
         {P(W_2\,|\,W_1,B_2,I_0)\cdot P(W_1\,|\,B_2,I_0)}
      \times  
\frac{P(B_1\,|\,I_0)}{P(B_2\,|\,I_0)} \label{eq:seq:2} \\
&=& \frac{P(W_2\,|\,B_1,I_0)}{P(W_2\,|\,B_2,I_0)}
    \times \frac{P(W_1\,|\,B_1,I_0)}{P(W_1\,|\,B_2,I_0)}
    \times \frac{P(B_1\,|\,I_0)}{P(B_2\,|\,I_0)} \label{eq:seq:3} \\
&=& \frac{P(W_2\,|\,B_1,I_0)}{P(W_2\,|\,B_2,I_0)}
    \times \frac{P(B_1\,|\,W_1,I_0)}{P(B_2\,|\,W_1,I_0)}\,, \label{eq:seq:4}
\end{eqnarray}
that, 
using the compact notation introduced above,
we can rewrite in the following enlighting forms. 
The first is 
[Eq.~(\ref{eq:seq:4})]
\begin{eqnarray}
O_{1,2}(W_1,W_2,I) &=& \tilde O_{1,2}(W_2,I) \times O_{1,2}(W_1,I)\,,
\end{eqnarray}
that is,
{\it the final odds after the first inference become the
initial odds of the second inference}
(and so on, if there are several pieces of evidence).
Therefore, beginning
from a situation in which $B_1$ was thirteen times more credible
than $B_2$ is exactly equivalent to having
started from unitary odds updated 
by a factor 13 due to a piece of evidence. 

The second form comes from Eq.~(\ref{eq:seq:3}): 
\begin{eqnarray}
O_{1,2}(W_1,W_2,I) &=& \tilde O_{1,2}(W_1,I) \times \tilde O_{1,2}(W_2,I) 
                       \times O_{1,2}(I) \label{eq:seq:3a}\\
                  &=&\tilde O_{1,2}(W_1,W_2,I) \times O_{1,2}(I)\,,
\label{eq:bf_mult} 
\end{eqnarray}
i.e.\footnote{Eq.~(\ref{eq:bf_mult}) follows from Eq.~(\ref{eq:seq:3a})
because a Bayes factor can be defined as the ratio 
of final odds over the initial odds, depending on the evidence.
Therefore
\begin{eqnarray*}
\tilde O_{1,2}(W_1,W_2,I) &=& \frac{O_{1,2}(W_1,W_2,I)}{ O_{1,2}(I)} 
= \tilde O_{1,2}(W_1,I) \times  \tilde O_{1,2}(W_2,I)\,.
\end{eqnarray*}
} 
\begin{eqnarray}
\tilde O_{1,2}(W_1,W_2,I) &=&  \tilde O_{1,2}(W_1,I) \times 
                              \tilde O_{1,2}(W_2,I)\,: 
\end{eqnarray}
{\it Bayes factors} due to 
{\it independent}\footnote{Probabilistic, or `stochastic', 
independence of the observations 
is related to the 
validity of the relation $P(W_2\,|\,W_1,B_i,I)=P(W_2\,|\,B_i,I)$,
that we have used above to turn 
Eq.~(\ref{eq:seq:2}) into Eq.~(\ref{eq:seq:3}) and that can be expressed,
in general terms as
\begin{eqnarray*}
P(E_2\,|\,E_1,H_i,I)=P(E_2\,|\,H_i,I)\,,
\end{eqnarray*}
i.e., {\it under the condition of a well precise
hypothesis} ($H_i$), the probability of the effect $E_2$
does not depend on the knowledge of whether $E_1$ 
has occurred or not. Note that, in general, 
although $E_1$ and $E_2$ are independent given $H_i$
(they are said to be {\it conditionally independent}),
they might be otherwise {\it dependent}, i.e.
$P(E_2\,|\,E_1,I_0)\ne P(E_2\,|\,I_0)$. 
(Going to the
example of the boxes, it is rather easy to grasp,
although I cannot enter in details here,
that, if we do not know the kind of box, the 
observation of $W_1$ changes our opinion about the box
composition and, as a consequence, 
the probability of $W_2$ -- see the examples in Appendix J)
}
pieces of evidence {\it multiply}. That is, 
two independent pieces of evidence ($W_1$ and $W_2$) are equivalent
to a single piece of evidence (`$W_1\cap W_2$'), whose Bayes factor
is the product of the individual ones. 
In our case $\tilde O_{1,2}(W_1\cap W_2,I)=169$.

In general, if we have several hypotheses $H_i$ and 
several {\it independent} 
pieces of evidence, $E_1$, $E_2$, \ldots, $E_n$, 
indicated all together as $\mvec{E}$,
then
Eq.~(\ref{eq:Bayes_factor}) becomes
\begin{eqnarray}
O_{i,j}(\mvec{E},I) &=& 
\left[\,\prod_{k=1}^n
\tilde O_{i,j}(E_k,I)
\right]  
\times O_{i,j}(I)\,,  \label{eq:product_Odds}
\end{eqnarray}
i.e. 
\begin{eqnarray}
\tilde O_{i,j}(\mvec{E},I) &=& \prod_{k=1}^n
\tilde O_{i,j}(E_k,I)\,,  \label{eq:product_Odds_1}
\end{eqnarray}
where $\prod$ stand for `product'
(analogous to $\sum$ for sums).

\subsection{How the independent arguments sum up in our judgement -- 
logarithmic updating and its interpretation}\label{ss:JL}
The remark that Bayes factors due to independent pieces
of evidence multiply together and the overall factor finally 
multiplies the initial odds suggests a change
of variables in order to play with additive quantities.\footnote{The 
idea of transforming a multiplicative updating into
an additive one via the use of logarithms 
is quite natural and seems to have 
been firstly used in 1878 by Charles Sanders Peirce~\cite{Peirce}
and finally introduced in the statistical practice mainly due
to the work of I.J. Good~\cite{Good}.
For more details see the Appendix E.}
This can be done taking the logarithm of both sides of 
Eq.~(\ref{eq:product_Odds}),
that then become
\begin{eqnarray}
\log_{10}[O_{i,j}(\mvec{E},I)] &=& 
\sum_{k=1}^n \log_{10}[\tilde O_{i,j}(E_k,I)] + \log_{10}[O_{i,j}(I)]\,,
\label{eq:sum_BF}
\end{eqnarray}
respectively, where
the base 10 is chosen for practical convenience because, 
as we shall discuss later,
what substantially matters are powers of ten of the odds.%

Introducing the new symbol JL, we can rewrite Eq.~(\ref{eq:sum_BF}) as
\begin{eqnarray}
\mbox{JL}_{i,j}(\mvec{E},I) &=&  \mbox{JL}_{i,j}(I) + 
\sum_{k=1}^n \Delta\mbox{JL}_{i,j}(E_k,I) 
\label{eq:Sum_Delta_JL}\\
&=&  \mbox{JL}_{i,j}(I) + 
 \Delta\mbox{JL}_{i,j}(\mvec{E},I) 
\label{eq:Sum_Delta_JL}
\end{eqnarray}
or 
\begin{eqnarray}
 \Delta\mbox{JL}_{i,j}(\mvec{E},I) &=& 
\mbox{JL}_{i,j}(\mvec{E},I) - \mbox{JL}_{i,j}(I)\,,
\label{eq:Delta_JL_diffJL}
\end{eqnarray}
where 
\begin{eqnarray}
\mbox{JL}_{i,j}(\mvec{E},I) &=& 
\log_{10}\left[O_{ij}(\mvec{E},I)\right] \label{eq:JL_allE}\\
\mbox{JL}_{i,j}(I) &=& 
\log_{10}\left[O_{i,j}(I) \right]\label{eq:JL_0}\\
\Delta\mbox{JL}_{i,j}(E_k,I) &=& 
\log_{10}\left[\tilde O_{i,j}(E_k,I)\right] 
 \label{eq:DeltaJL_k} \\
 \Delta\mbox{JL}_{i,j}(\mvec{E},I)  &=& 
\sum_{k=1}^n \Delta\mbox{JL}_{i,j}(E_k,I)\,.
  \label{eq:Sum_DeltaJL_k}
\end{eqnarray}
\begin{figure}
\centering\epsfig{file=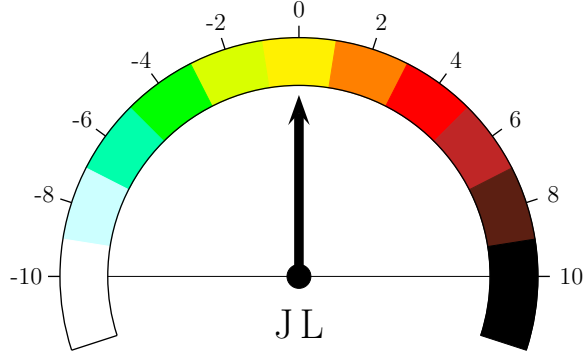,clip=,width=0.5\linewidth}
\caption{\small \sf Judgement leaning.}
\label{fig:jl}
\end{figure}
The letter `L' in the symbol is to remind 
{\it logarithm}. But it has also 
the mnemonic meaning of 
{\it leaning}, in the sense of `inclination' or `propension'. 
The `J' is for {\it judgment}. Therefore
`JL' stands for {\it judgement leaning},
that is an inclination of the judgement,
an expression I have taken the liberty to introduce,
using words not already engaged in probability and statistics,
because in these fields many controversies are due to different
meanings attributed to the same word, or expression,
by different people (see Appendices B and G 
for further comments). 
JL can then be visualized as  the indicator of the 
`justice balance'\footnote{I have realized only later that  
JL sounds a bit like `jail'. That might be not so bad, 
if $H_1$ to which $\mbox{JL}_{1,2}(E_k)$ refers stands for `guilty'.}
(figure \ref{fig:jl}),
that displays zero if there is no unbalance, but it could move
to the positive or the negative side depending on the 
weight of the several arguments pro and con.
The role of the evidence is to vary the JL indicator by quantities
$\Delta\mbox{JL}$'s equal to base 10 logarithms of the
Bayes factors, that have then a meaning of 
{\it weight of evidence}, an expression due to 
Charles Sanders Peirce~\cite{Peirce} (see Appendix E).

But the judgement is rarely initially unbalanced.
This the role of $\mbox{JL}_{i,j}(I)$, that can be considered 
as a a kind of {\it initial weight of evidence} due to our prior knowledge
about the hypotheses $H_i$ and $H_j$ [and that could even
be written as  $\Delta\mbox{JL}_{i,j}(E_0,I)$, 
to stress that it is related
to a 0-th piece of evidence]

To understand the rationale behind a possible uniform treatment
of the prior as it would be a piece of evidence, 
let us start from a case in which you now 
{\it absolutely nothing}. For example \underline{You} have to state
your beliefs on which  of my friends, Dino or Paolo, 
will first run next Rome marathon. It is absolutely reasonable
you assign to the two hypotheses  equal probabilities, i.e. $O_{1,2}=1$, 
or $\mbox{JL}_{1,2}=0$ (your judgement is
perfectly balanced). This is because in \underline{Your}
brain these
names are only possibly related to Italian males. Nothing
more. (But nowadays search engines over the web allow to
modify your opinion in minutes.)

As soon as you deal with {\it real} hypotheses of your interest,
things get quite different.
It is in fact very rare the case in which the hypotheses
tell you not more than their names.
It is enough you think at the hypotheses  `rain' or `not rain', 
the day after you read these lines in the place where you live. 
In general the information you have in your brain
related to the hypotheses of your interest can be considered 
the initial piece of evidence {\it you} have, 
usually different from that somebody else
might have
(this the role of $I$ in all our expressions). 
It follows that prior odds of 10 ($\mbox{JL}=1$) will influence your 
leaning
towards one hypothesis, exactly like 
unitary odds ($\mbox{JL}=0$) followed by a Bayes factor of 
10 ($\Delta \mbox{JL}=1$).
This the reason they enter on equal foot when 
``balancing arguments''
(to use an expression \`a la Peirce -- see the Appendix E)
 pro and against hypotheses.
\begin{table}[t]
\begin{center}
\begin{tabular}{|c|c|c|c|c|c|c|}
\hline
Judg. leaning  & Odds(1:2)  & $P(H_1)$  & &  
Judg. leaning  & Odds(1:2)   & $P(H_1)$ \\ 
\,[$\mbox{JL}_{1,2}$]    &   [$O_{1,2}$]& (\%) & & $\mbox{JL}_{1,2}$] &   [$O_{1,2}$]& (\%)   \\
\hline
$0$ & 1.0     & 50   &&  &    &           \\
$-0.1$ & 0.79  & 44  &&  $0.1$ & 1.3   &  56  \\   
$-0.2$ & 0.63  & 39  &&  $0.2$ & 1.6   &  61  \\ 
$-0.3$ & 0.50  & 33  &&  $0.2$ & 2.0   &  67  \\ 
$-0.4$ & 0.40  & 28  &&  $0.4$ & 2.5   &  71  \\  
$-0.5$ & 0.32  & 24  &&  $0.5$ & 3.2   &  76  \\ 
$-0.6$ & 0.25  & 20  &&  $0.6$ & 4.0   &  80  \\ 
$-0.7$ & 0.20  & 17  &&  $0.7$ & 5.0   &  83  \\
$-0.8$ & 0.16  & 14   && $0.8$ & 6.3   &  86  \\
$-0.9$ & 0.13  & 11   && $0.9$ & 7.9   &  89  \\
$-1.0$ & 0.10  & 9.1  && $1.0$ & 10    &  91  \\
$-1.1$ & 0.079  & 7.4  && $1.1$ & 13   &  92.6  \\
$-1.2$ & 0.063 & 5.9  && $1.2$ & 16    &  94.1    \\ 
$-1.3$ & 0.050  & 4.7  && $1.0$ & 20    &  95.2  \\
$-1.4$ & 0.040 & 3.8  && $1.4$ & 25    &  96.2    \\ 
$-1.5$ & 0.032 & 3.1  && $1.5$ & 32    &  96.9    \\ 
$-1.6$ & 0.025 & 2.5  &&  $1.6$ & 40    &  97.5  \\
$-1.7$ & 0.020 & 2.0  &&  $1.7$ & 50    &  98.0  \\
$-1.8$ & 0.016 & 1.6  &&  $1.8$ & 63    &  98.4  \\
$-1.9$ & 0.013 & 1.2  &&  $1.9$ & 80    &  98.8  \\
$-2.0$ & 0.010 & 1.0  &&  $2.0$ & 100   &  99.0  \\
\hline
\end{tabular}
\caption{\small \sf A comparison between 
probability, odds and judgement leanings 
%[write javascript to practice!]
}
\label{tab:JL_BF_P}
\end{center}
\end{table}

Finally, table \ref{tab:JL_BF_P} compares judgements leanings, odds and
probabilities, to show that the human sensitivity to belief 
(that is something like 
Peirce's {\it intensity of belief} -- see Appendix E) 
is not linear with probability. 
For example, if we assign
probabilities of 44\%, 50\% or 56\% to events $E_1$, $E_2$ and $E_3$
we do not expect one of them
really more strongly than the others, 
in the sense that we  are not much 
{\it surprised} of any of the three occurs. But the same differences
in probability produce quite different {\it sentiment of surprise} 
if we shift the probability scale (if they were, instead, 1\%, 7\% and 13\%,
we would be highly surprised if $E_1$ occurs). 

Similarly 99.9\% probability on $H$ is substantially different from 99.0\%,
although the difference in probability is `only' 0.9\%. This is well 
understood, and in fact it is known that the best way to express 
the perception of probability values very close to 1 is to think to the opposite
hypothesis $\overline H$, that is 0.1\% probable in the first case and 1\%
probable in the second -- we could be quite differently surprised 
if $H$ does not result to be true in the two cases!\footnote{The
`switch of perspective' from $E$ to $\overline H$ is done 
in a way somewhat automatic if, instead of the probability, we take
the logarithm of the odds, for example our JL (obviously the 
base of the logarithm is irrelevant). Since 
$\mbox{JL}_H(I)=\log_{10}[P(H\,|\,I)/P(\overline H\,|\,I)]$, 
in the limit $P(H\,|\,I)\rightarrow 0$ we have that
$\mbox{JL}_H(I)\approx \log_{10}[P(H\,|\,I)]$, while 
 the limit $P(H\,|\,I)\rightarrow 1$ it is 
$\mbox{JL}_H(I)\approx - \log_{10}[P(\overline H\,|\,I)]$.
\label{fn:logP}
}

From the table we can see that the {\it human resolution} is about 1/10
of the JL, 
although this does not imply that a probability value of
53.85\% ($\mbox{JL}=0.0670$) cannot be stated. It all depends
how this value has been evaluated and what is the purpose of 
it.\footnote{This is more or less what happens in measurements. 
Take for example the probabilities that appears in the 
$E_1$ `monitor' of figure \ref{fig:bn_mon_0}: 53.85\%
for white and 46.15\% for black. This is like to say that
two bodies weigh 53.85\,g and 46.15\,g, as resulting
from a measurement with a precise balance (the Bayesian network 
tool described in Appendix J applied to the box toy model 
is the analogue of the precise balance). 
For some purposes two, three and even four significant digits 
can be important. But, anyhow, as far as
our perception is concerned, not only the least digits are
absolutely irrelevant but we can hardly distinguish between 
54\,g and 46\,g.
}

\subsection{Recap of the section}
This section had the purpose of introducing the so-called
Bayesian reasoning (that is, in reality, nothing more than just
probabilistic reasoning) with an aseptic, simple
example, that  shows however the ingredients needed
to update our opinion on the light of new
observations.
At this point the role of the priors and 
of the evidence in forming our
opinion about the hypotheses of interest should be  clear. 
Note also how I have used on purpose 
several expressions to mean essentially the same thing, 
expressions that involve words such as `probability', `belief',
`plausibility', `credibility', `confidence', and so on. 

%\newpage
\section{Weight of priors and weight of evidence in real life}
\label{sec:real_life}
The box example used to introduce the Bayesian 
reasoning was particularly simple for two reasons. 
First, the updating factor was calculated 
from elementary probability rules in an `objective
way' (in the sense that everybody would agree on 
a Bayes factor of 13, corresponding to a $\Delta$JL of 1.1). 
Second, also the
prior odds $n_1/n_2$ were univocally determined by the 
formulation of the problem. 

In real life the situations are never so simple.
Not only priors can differ a lot from a person
to another. Also the probabilities that enter
the Bayes factor might not be the same for everybody.
Simply because they are probabilities, and
probabilities, meant as degree of belief,
have an intrinsic {\it subjective}
nature~\cite{BdF}. 
The very reason for this trivial remark 
(although not accepted by everybody, because of 
ideological reasons) is that probability depends
on the available information and
-- fortunately! -- there are no two identical brains
in the world, made exactly the same way and
 sharing exactly the same information. 
Therefore, the same event is not expected with the same 
security by different subjects, and the same hypothesis
is not considered equally credible.\footnote{The following 
quotes can be rather enlighting,
especially for those who think they think, 
just for educational reasons, `they have to be frequentist':
\begin{quote}
{\sl ``Given the state of our knowledge about everything that could 
possibly have any bearing on the coming true of a certain event
(thus {\it in dubio}: of the sum total of our knowledge), 
the numerical probability $p$ of this event is to be a real number
by the indication of which we try in some cases to set up a quantitative
measure of the strength of our conjecture or anticipation, founded 
on the said knowledge, that the event comes true. \\
\ldots \\
Since the knowledge may be different with different persons or with
the same person at different times, they may anticipate the same event
with more or less confidence, and thus different numerical probabilities
may be attached to the same event.\ \ldots\  Thus 
whenever we speak loosely of the `probability of an event,'
it is always to be understood: probability with regard to a certain
given state of knowledge.''}~\cite{SchroedingerA}
\end{quote}
} 

At most degrees of belief
can be {\it inter-subjective}, because in 
many cases there are people or entire communities 
that share the same initial beliefs (the same {\it culture}),
reason more or less
the same way (similar {\it brains} and similar {\it education})
and have access to the {\it same data}.
Finally, there are stereotyped
`games' in which probabilities can even be {\it objective},
in the sense that everybody will agree on its value.
But these situations have to be considered the exceptions rather than the rule
(and even when we state with great security that the probability
of head tossing a regular coin is exactly 1/2, we forget 
it could remain vertically, a possibility usually
excluded but that I have personally experienced 
a couple of times in my life.)

Therefore, although educational games with boxes and balls 
might be useful to learn the grammar and syntax
of probabilistic reasoning, at a given point we need to
move to real situations.

\subsection{Assessing subjective degrees of beliefs -- virtual bets}
A good way to force experts to provide 
the initial beliefs they
have formed in their minds, 
elaborated somehow by their `educated intuition'
(see Appendix C), is to propose them a virtual lottery, 
in which they can choose the event on which
to bet to win a rich prize. 
One is the event of interest (let us call it $A$), 
the other one is a simpler one, based on  
coins, urns, dice or card games. The latter can be considered
a kind of `standard',
or a `reference' (as it is
done in measurements to calibrate instruments),
whose probability is the same 
for everyone.
We can ask ourselves (or the experts), for example, if
we (or they) prefer
to bet on $A$ rather than on head resulting from a regular coin; 
or 
on white extracting a ball
from a box containing 100 balls, 90 of which white;
and so on. 

Obviously, none can state initial odds with very
high precision.\footnote{Those 
who are not familiar with this
approach have understandable
initial difficulties and risk to be at lost. 
A formula, they might argue, 
can be of practical use only if we can replace 
the symbols by numbers, and in pure mathematics a number 
 is a well defined
object, being, for example, 49.999999 different from 
50.
Therefore, they might conclude that, being unable to 
choose {\it the} number, the above formulae, that 
seem to work nicely in die/coin/ball games, are
useless in other domains of applications (the most
interesting of all, as it was clear already  centuries ago
to Leibniz and Hume). 
But in the realm of uncertainty things go quite differently,
as everybody understands, apart from hypothetical 
Pythagorean monks living in a ivory monastery. 
For practical purposes not only 49.999999\% is `identical'
to 50\%, but also 49\% and 51\% 
give to our mind essentially the same 
expectations of what it could occur.
In practice we are interested to 
understand if somebody else's degrees of belief
are low, very low, high, very very high, ad so on.
And the same is what other people expect from us.
} 
But this does not matter (table \ref{tab:JL_BF_P} 
can help to get the point). 
We want to understand
if they are of the order of 1 (equally likely), 
of the order of a few units (one is a bit more likely than the other one), 
or  of suitable powers of 10 (much more or much
less likely than the other one). If one has doubts about the final
result, one can make a `sensitivity analysis', i.e. vary the 
value inside a wide but still  believable range and check
how the result changes. The sensitivity 
(or insensitivity) will depend also on the 
other pieces of evidence to draw the final conclusion.
Take for example two different evidences, characterized by 
Bayes factors of $H_1$ versus $H_2$ 
very high (e.g. $10^4$)
or very small (e.g. $10^{-4}$),  
corresponding to $\Delta$JL's of $+4$ or $-4$,  respectively
(for the moment we assume all subjects agree 
on the evaluation of these factors). 
Given these values, it is easy to check that,
for many practical purposes,
the conclusions will be the same 
even if the initial odds are in the range $1/10$ to $10$,
i.e. a JL between $-1$ and $+1$, that can be stated
as $\mbox{JL}_{1,2}(E_0)=0\pm 1$. Adding 
`weights of evidence' of $+4$ or $-4$, we get
final JL's of $4\pm 1$ or $-4\pm 1$, 
respectively.\footnote{
That is, the final probability of $H_1$ would range between 99.90\%
and 99.999\% in the first case, 
between 0.001\% and 0.1\% in the second one,
making us `practically sure' of either hypothesis 
in the two cases.}

The limit case in which the Bayes factor
is zero or infinity (i.e. $\Delta$JL's $-\infty$ or $+\infty$)
makes the conclusion absolutely
independent from priors, as it seems obvious.

\subsection{Beliefs versus frequencies}\label{ss:belief_freq}
At this point a remark on the important (and often misunderstood)
issue of the relation between degrees of beliefs and relative 
frequencies is in order.

The concept of subjective probability does not
preclude the use of
{\it relative frequencies} in the reasonings. 
In particular, beliefs can be evaluated from the relative 
frequencies of other events, analogous to the 
one of interest, have occurred in the past. 
This can be done roughly (see Hume's quote in Appendix B)
or in a rigorous way, 
using probability theory 
under well defined assumptions
(Bayes' theorem applied
to the inference of 
the parameter $p$ of a binomial distribution).

Similarly, if we believe that a given event will
occur with 80\% probability, it is absolutely correct
to say that, if we think at a large number of analogous 
independent events that we consider equally probable,
we expect that in about 80\% of the cases the event
will occur. This also comes from probability theory
(Bernoulli theorem). 

This means that, contrary to what one reads often 
in the literature and on the web, 
evaluating probabilities from past frequencies
and expressing beliefs by expected frequencies does
not imply to adhere to {\it frequentism}.\footnote{Sometimes
frequency is even confused with `proportion'
when it is said, for example, that the probability
is evaluated thinking how many persons in a given
population would behave in a given way, or have a 
well defined character.}  The importance
of this remark in the context of this paper is that 
people might find natural, for their reasons, 
to evaluate and to express beliefs this way,
although they are perfectly aware that the event
about they are reasoning is unique. For further comments see
Appendix B.

\subsection{Subjective evaluation of Bayes factors}
As we have mentioned above, and as we shall see later,
not always the evaluation of updating factors can be done
 with the help of mathematical formulae like in the box example.
However, we can make use of the virtual bet in this case too,
remembering that a Bayes factor can be considered 
as the odds due a single piece of evidence,
provided the two hypotheses are considered otherwise
equally likely (hence, let us remember, the symbol $\tilde O$ used here
to indicate Bayes factors).

\subsection{Combining uncertain priors and uncertain weights of evidence}
\label{ss:combining}
When we have set up our problem, listed the 
pieces of evidence pro and con, including the $0$-th one
(the prior), and attributed to each of them a weight of 
evidence, quantified by the corresponding $\Delta$JL's,
we can finally sum up all contributions.

\begin{figure}
\centering\epsfig{file=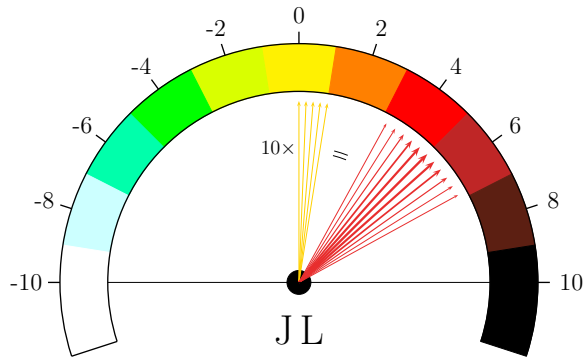,clip=,width=0.5\linewidth}
\caption{\small \sf Combined effect of
10 weak and `vague' pieces of evidence, 
each yielding a $\Delta$JL of $0.5\pm 0.5$ (see text).}
\label{fig:jl_1}
\end{figure}

As it is easy to understand, if the number of pieces of evidence
becomes large, the final judgment can be rather precise 
and far from being perfectly balanced, 
even if each contribution is weak and even uncertain.
This is an effect of the famous 
`central limit theorem' that dumps
the weight of the values far from the 
average.\footnote{The reason behind it is rather easy 
to grasp. When we have uncertain beliefs
it is like if {\it our mind oscillates} among possible
values, without being able to choose an exact value.
Exactly as it happens when we try to guess, just by eye,
the length of a stick, the weight of an object or a 
temperature in a room: extreme values are 
promptly rejected, and our judgement oscillates in an interval,
whose width depends on our estimation ability, based on 
previous experience. 
Our guess will be somehow the center of the interval.
The following minimalist example helps to understand 
the rule of combination of uncertain evaluations.
Imagine that the (not better
defined) quantities $x$  and $y$ might each 
have, in our opinion,
the values 1, 2 or 3, 
among which we are unable to choose.
If we now think of a $z=x+y$, its value can then range between 
2 and 6. But, if our mind oscillates uniformly and independently
over the three possibilities of $x$ and $y$, the oscillation
over the values of $z$ is not uniform. The reason is that $z=2$
can is only related to $x=1$ and $y=1$. Instead, we think at 
$z=3$ if we think at $x=1$ and $y=2$, or at $x=2$ and $y=1$. 
Playing with a cross table of possibilities, it is rather
easy to prove that $z=4$ gets
a weight three times larger than that of $z=2$.
We can add a third quantity $v$, similar to $x$ and $y$,
and continue the exercise, understanding  then the essence of 
what is called in probability theory {\it central limit theorem},
which then applies also to the weight of our JL's.
[Solution and comment: if $w=z+v$, the weights of the 7 possibilities, from 3 to 9
are in the following proportions: 1:3:6:7:6:3:1. Note that, contrary
to $z$, the weights do not go linearly up and down, but there is
a non-linear concentration at the center. When many variables of this
kind are combined together, then the distribution
of weights exhibits the well known bell shape of the Gaussian distribution.
The widths of the red arrows in figure \ref{fig:jl_1}
tail off from the central one according to a Gaussian function.]
\label{note:mind_oscillations}
} 
Take for example the case of 10 JL's each 
uniformly\footnote{It easy to understand that if the 
judgement would be uniform in the odds, ranging then from 
1 to 10, the conclusion could be different. Here it
is assumed that the `intensity of belief'\cite{Peirce} 
is proportional to the logarithm of the odds, as 
extensively discussed in Appendix E.}
ranging between 0 and 1, i.e. 
$\Delta\mbox{JL}_{1,2}(E_i,I)=0.5\pm 0.5$.
Each piece of evidence is marginal, but the sum leads to
a combined  $\Delta\mbox{JL}_{1,2}(\mvec{E},I)$ of $5.0\pm 1.8$, 
where ``$[3.2,6.8]$'' defines now an effective range 
of leanings\footnote{Using the language of footnote
\ref{note:mind_oscillations}, this is the range in which 
the minds oscillate in 95\% of the times when thinking of 
$\Delta\mbox{JL}_{1,2}(\mvec{E},I)$.},
as depicted in figure \ref{fig:jl_1}. 
Note that in this graphical representation
the 5 yellow arrows (the lighter ones if you are reading
the text in black/white) 
do not represent individual values of JL,
but its interval. These arrows have all the same width 
to indicate that the exact value is indifferent to us. 
The red arrow have instead different widths to indicate that
we prefer the values around 5 and the preference goes down as we 
move far form it. The 12 arrows only indicate an effective range,
because the full range goes from 0 to 10, although $\Delta$JL values
very far from 5 must have negligible weight in our reasoning.

\subsection{Agatha Christie's
``three pieces of evidence''}\label{ss:Agatha}
As we have seen, 
a single evidence, 
yielding a Bayes factor of the order of 10, or a $\Delta$JL around 1,
is not a strong evidence. But many individual, independent
pieces of evidence of that weight 
should have much a greater consideration
in our judgement.

This is, somehow, 
the rational behind Agatha Christie's
``three pieces of evidence''.
However it is worth remarking that 
something is to
say there is a rational behind this expression, that
can be used as a rough rule of thumb, something else
is to take it as a `principle', as it is often supposed
in the Italian dictum ``tre indizi fanno una prova''.
First, 
pieces of evidence are usually not `equally strong', in the sense
they do not carry the same weight of evidence
and sometimes even several pieces of evidence are 
not enough.\footnote{I wish 
 judges state Bayes factors 
of {\it each} piece of evidence, as vaguely as 
they like 
(much better than telling nothing! -- Bruno
de Finetti was used to say that ``{\it  it is better
to build on sand that on void}''),
instead of saying that somebody is 
guilty ``behind any reasonable doubt'' --
and I am really curious to check to what 
degree of belief that 
level of doubt corresponds!}
Second, the prior -- that is our `0-th evidence -- 
can completely balance the weight of evidence. Finally,
we have also to remember that sometimes
they are not even completely independent, in which case the
product rule is not any longer 
valid.\footnote{What to do in this case? 
As it easy to imagine, when the structure of dependencies
among evidences is complex, things might become 
quite complicated. Anyway, if one is able to isolate
two o more pieces of evidence that are 
correlated with themselves 
(let they be $E_1$ and $E_2$), then, one can consider
the joint event $E_{1\&2}=E_1\cap E_2$ 
as the effective evidence to be used.
In the extreme case in which $E_1$ implies logically $E_2$
(think at the events `even' and '2' rolling a die),
then $P(E_2\,|\,E_1,I)=1$, from which it follows that
$P(E_1\cap E_2\,|\,I)=P(E_1\,|\,I)$: the second evidence 
$E_2$ is therefore simply superfluous.
}

A final remark on the combination of pieces of
evidence is still in order.
From a mathematical point of view there is no difference
between a single piece of evidence 
yielding a tremendous Bayes factor
of $10^{10}$ ($\Delta\mbox{JL}=10$) and ten independent pieces of evidence,
each having the more modest Bayes factor of 10 
($\Delta\mbox{JL}=1$). 
However, I have somehow the
impression (mainly got from media and from fiction, since I have no
direct experience of courts) 
that the first is considered as {\it the} incriminating evidence
(the `smoking gun'),
while the ten weak pieces of evidence are just taken as some
{\it  marginal  indications},
that all together are not as relevant as the single incriminating
`proof'. Not only this reasoning is mathematically incorrect, 
as we have learned, but, if I were called to state my opinion 
on the two sets of evidence, I had no doubt to consider 
the ten weak pieces
of evidence more incriminating than the single `strong' one, 
although they seem to be formally equivalent. Where is the point?
In all reasonings done until now we have focused on the 
weight of evidence, assuming each evidence is a true
and not a fake one, for instance incorrectly reported, 
or even fabricated by the investigators.
In real cases one has to take
into account also this 
possibility.\footnote{When we are called to make
critical decisions 
even very remote hypotheses, 
although with very low probability,
should be present to our minds  
-- that is 
Dennis Lindley's 
{\it Cromwell's rule}~\cite{Lindley}. 
[The very recent news from New York 
offer material for reflection~\cite{AdilPolanco}.] 
\label{fn:Cromwell}}
As a consequence, if there 
is any slight doubt on the validity of each piece of evidence,
it is rather simple to understand that the single evidence is
somewhat weaker than the ten ones all together
(Agatha Christie's three pieces of evidence are in qualitative agreement
with this remark).
For further details see Appendix I.

\subsection{Critical values for guilt/innocence -- 
Assessing beliefs versus making decisions}
At this point a natural question raises spontaneously.
What is the 
possible threshold of odds or of JL's to 
condemn or to absolve somebody? 
This is a problem of a different kind. It is not just 
a question of believing something, but on deciding which action
to take. 

Decision issues are a bit more complicate than probability
ones. Not only they inherit all probabilistic questions, 
but they need careful considerations of all possible benefits
and losses resulting from the action. I am not a judge
and fortunately I have never been called to join 
a popular jury, on the validity of which I have, by 
the way, quite some doubts.\footnote{Again, my impression comes 
from media, literature and fiction, but I cannot see how
`casual judges' can be better than professional ones
to evaluate all elements of a complex trial, or how to
distinguish sound arguments from pure rhetoric of the 
lawyers. This is particularly true when the `network of evidences'
is so intricate that even well trained human minds
might have difficulties, and {\it artificial intelligence}
tools would be more appropriated (see Appendices C and J).}
So I do not know exactly how they make their decisions, 
but personally, being 99\% confident that somebody
is guilty (that is
a JL of 2), I would
not behave the same way if the person is accused of a 
`simple' crime of passion, or of being a Mafia or a serial killer.

\subsubsection{I have a dream}
I hope judges know  what they do,
but I wish one day they will finally state somehow, in a quantitative
way, with all possible uncertainties, 
the beliefs they have in their mind, the individual 
contributions they have considered and the 
society benefits and losses taken into account
to behave the way they did.

\section{Columbo's priors versus jury's priors}\label{sec:Columbo_priors}
Going back to Columbo's episode, the prior of interest here
is the probability
that Peter Galesco killed his wife, taking 
into account `all' pieces of evidence {\it but} those
deriving from the last scene.

It is interesting to observe how
{\it probabilities} change as the story goes on. 
Different characters
develop different opinions, depending on their 
previous experience, on the information they get
and on their capability to process the information quickly.
Also each spectator
forms his/her own opinion, although all of them get
virtually the same `external' pieces of information
(that however are combined  with internal pre-existing 
ones, whose combination and rapidity of combination 
depend on many other internal things 
and environmental conditions)
-- and this is part of the fun of watching 
a thriller with friends.

\subsection{Columbo's priors}
By definition a person suspected by a detective
is not just anybody, whose name was extracted
at random from the list of citizens in the region
where the crime was committed. 
Police does not like to lose time, money and reputation,
if it does not have valid suspicions, and
investigations proceed in different directions, 
with priorities proportional to the chance
of success. 
The probabilities of the various hypotheses go up
and down as the story goes on, and
an alibi or a witness could drop a probability
to zero (but policemen are aware of  fake alibis
or lying witnesses).

If we see Columbo loosing sleep following some hints,
we understand he has strong suspicions. Or, 
at least,
he is not convinced of the official version of the
facts, swallowed  instead by his colleagues: some 
elements of the puzzle do not fit nicely together
or, told in probabilistic terms, the 
{\it network of beliefs}\footnote{See Appendices C and J.}
he has in mind\footnote{Obviously, 
saying Columbo has a network of beliefs in his head,
I don't mean he is thinking at these mathematical tools.
On the other way around, these tools try to model
the way we reason, with the advantage they 
can  better handle complex
situations (see Appendices C and J).
\label{note:BN}
} 
makes him {\bf highly confident} 
that the suspected person is guilty.

\subsection{Court priors}
But a policeman is not the court that finally 
returns the verdict. Judges tend, by their experience, 
to trust policemen, but they cannot have 
exactly the same information the direct investigators 
have, that is not limited to what appears in the 
official reports.

Columbo might have formed
his opinion on instinctive reactions of Galesco,
on some photographer's hints of smile or on nervous replies
to fastidious questions, and so on, all little things
the lieutenant knows they cannot enter in the formal
suit.\footnote{There is, for example, the interesting 
case of the clochard who
was on the scene of the crime and, although still drunk, 
tells, among other verifiable things, to have heard two gun shots
with a remarkable
time gap in between, something in absolute contradiction
with Galesco reconstruction of the facts, in which he 
states to have killed Alvin Deschler, 
that he pretends to be the kidnapper 
and murderer of his wife, for self-defense,
thus shooting practically simultaneously with him. 
Unfortunately, days after, 
when the clochard is interviewed by Columbo, he says, 
apparently honestly,
to remember nothing of what happened the day of the crime, 
because
he was completely drunk. He confesses he
doesn't even remember what he 
declared to the police immediately after. 
Therefore he could never be able to testify
in a court. However, it is difficult an investigator 
would remove such a 
piece of evidence from his mind, a piece of evidence
that fits well with the alternative hypothesis 
that starts to account better for many other major and minor 
details.
He knows he cannot present it to the court, 
but it pushes him to go further, looking for more 
`presentable' pieces of evidence, and possibly 
for conclusive
proofs.}
We can form ourself 
an idea about the prior probability
that the court can assign to the hypothesis that 
the photographer is guilty from the reaction
of Columbo's colleagues and superiors, 
who try to convince him the case is settled: 
{\bf quite low}.

\subsection{Effect of a Bayes factor of 13}
To evaluate how a new piece of evidence modifies these
levels of confidence, we need to quantify somehow
the different priors. Since, as we have seen above, 
what really matters in these cases are the powers
of ten of the odds, we could place Columbo's ones
in the region $10^{2}$-$10^{3}$, the hypothetical
jury ones around $10^{-2}$, perhaps up to
$10^{-1}$. Multiplying these
values by 13 we see that, while the lieutenant 
would be practically sure Galesco is guilty,  
the jury component could hardly 
reach the level of a sound 
suspicion.

Using the expressions of subsection
\ref{ss:JL}, a Bayes factor of 13 corresponds to 
$\Delta\mbox{JL}=1.1$, that, added to initial leanings 
of $\approx 2.5\pm 0.5$ (Colombo) 
and $\approx -1.5 \pm 0.5$ (jury), could lead to
combined $\mbox{JL}$'s of $3.6\pm 0.5$ or 
 $-0.4 \pm 0.5$ in the two cases.

However, although such a small weight of evidence 
is not enough, by itself, to condemn a person,
I do not agree that
{\it ``that kind of evidence would 
never stand up in court''}\cite{NS} for the reasons
expounded in section~\ref{ss:combining}. 

Nevertheless, my main point in this paper is not that even 
such a modest piece of evidence 
{\it should} stand up in court (provided it is not the only one),
but rather that the weight of evidence 
provided by the rash Galesco's act is not 1.1, 
but much higher, infinitely higher.

\section{The weight of evidence of the
 full sequence of actions}\label{sec:finale}
In the previous section we have done the exercise
of assuming a Bayes factor
of 13, that is a weight of evidence of 1.1, 
as if taking that camera would be the same as
extracting a ball from a box,
as in the introductory example.
But does this look reasonable? 

\subsection{The `negative reaction'}
\begin{figure}
\hspace{-1.0cm}\epsfig{file=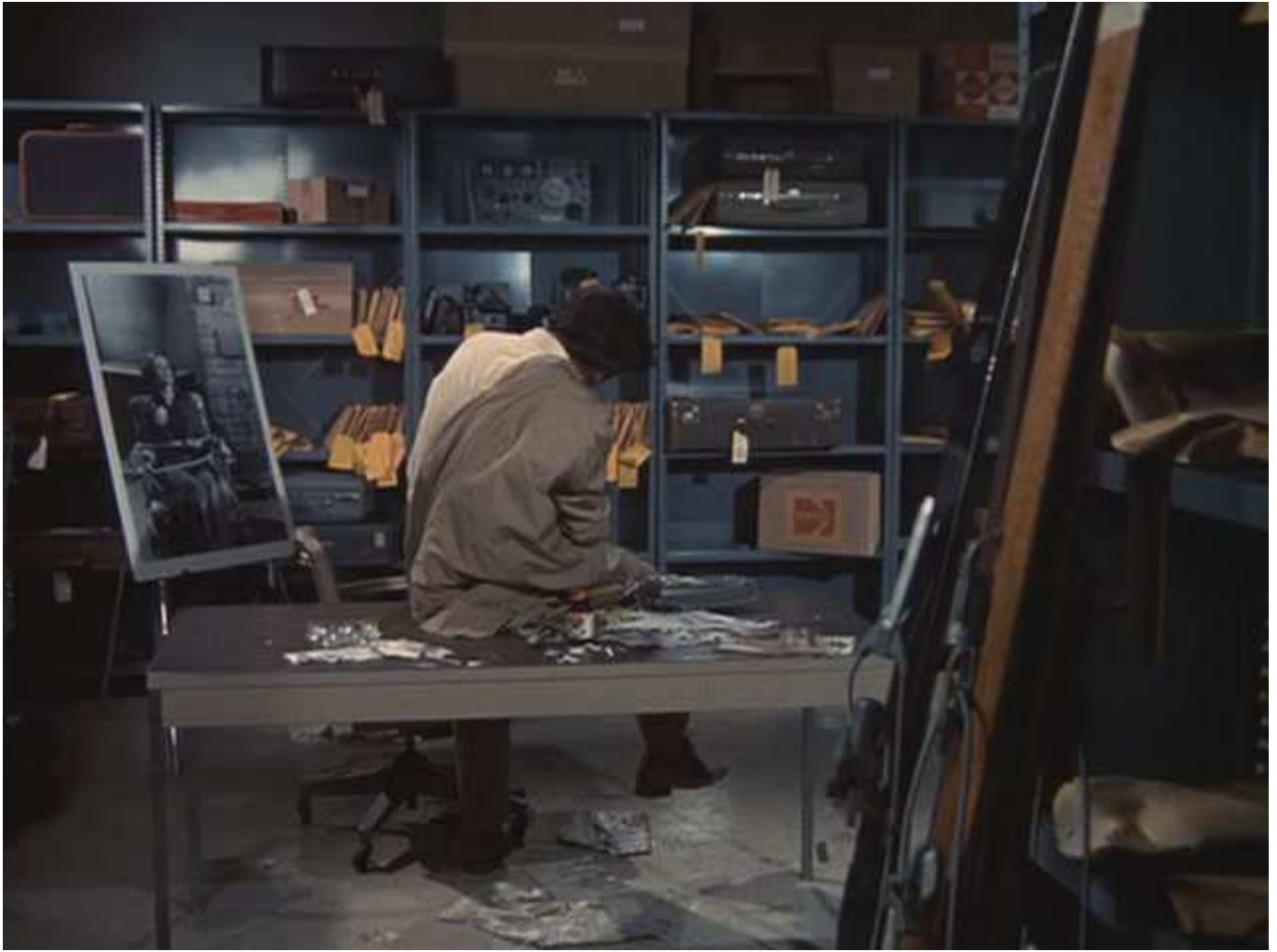,clip=,width=1.2\linewidth}
%http://raindogg.tistory.com/2856
\caption{\small \sf Final photogram of {\it Negative Reaction}. 
The incriminating camera is on the desk.
The remaining twelve are in the shelf just behind Columbo's head.
Desk and floor are full of the bits and pieces of the 
newspaper from which
the lieutenant tried to reproduce the kidnapper's note.}
\label{fig:last_photogram}
\end{figure}

Let us summarize what happens in the last scene of 
the episode.
\begin{itemize}
\item Galesco is suddenly taken to the police station, 
      where he is waited for by Columbo, who receives him not
      in his office but in a kind
      of repository containing shelfs full of objects
      (see figure \ref{fig:last_photogram}), including 
      the cameras in question, a few of which are 
      visible behind Columbo, although nobody mentions them.
\item Columbo starts arguing about the rest of the 
      newspaper found in Deschler's motel room 
      and used to cut out the words glued in the kidnapper's note.
      The missing bits and pieces 
      support the hypothesis that the collage was not done 
      by Deschler. Galesco, usually very prompt in suggesting
      explanations to Columbo's doubts and insinuations, 
      is surprised by the frontal attack of the lieutenant,
      who until that moment only expressed him a series of doubts.
      He gets then quite upset.
\item Immediately after, Columbo announces his final proof,
      meant to destroy Galesco's alibi. He has prepared 
      a giant enlargement  of the picture of  Mrs Galesco
      taken by the murderer 
      just before she was killed. The photograph shows clearly
      a clock indicating
      exactly ten (A.M.), time at which the lady had to be 
      with her husband, while Deschler had a very solid alibi,  that morning
      doing the driving test to get his licence. 
\item The expert photographer refuses this new reconstruction,
      on the ground that, he claims, there
      is a mistake in the 
      enlargement, in which, he says, the picture has been erroneously 
      reversed, thus 
      transforming the original  2:00 (P.M.) into 10:00
      (obviously, the analog clock had no digits,
      but just marks to indicate the hours). He asks 
      then Columbo to check on the original.
\item But Columbo acts very well in pretending he destroyed
      the original by accident when he was in the 
      dark room to supervise the work
      (his often goofy way 
      to behave makes the thing plausible).
\item This clever move is able to stress the otherwise always
      lucid Galesco, who suddenly thinks
      he is going to fall into a trap, based on a false,
      incriminating evidence fabricated by the police.
      He gets then so  
      nervous to loose control 
      and, with a kind of desperate jump of 
      a feline who sees itself lost, does his fatal mistake.
\begin{figure}
\centering\epsfig{file=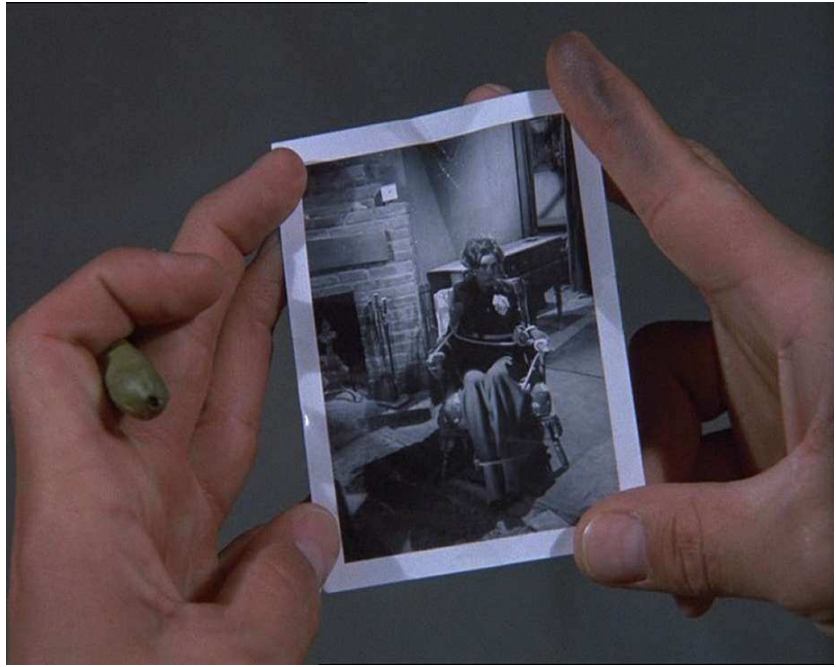,clip=,}
%http://raindogg.tistory.com/2856
\caption{\small \sf The discarded photograph in Columbo's hands.}
\label{fig:second_picture}
\end{figure}
\item Suddenly he has a kind of inspiration.
      He says the negative can 
      prove the picture has been reversed.\footnote{In reality
      he has several ways out, not depending on that negative
      (this could be a weak point of the story, but it is
      plausible, and the dramatic force of the action induces also
      TV watchers to neglect this particular, as my friends and I have
      experienced): 
      \begin{enumerate}
      \item 
      He knew 
      Columbo owns a second picture, discarded by the killer
      because of minor defects
      and left on the crime scene.
      (That was one of the several hints against Galesco, because only 
      a maniac photographer -- and certainly not 
      Alvin Deschler -- would care of 
      the artistic quality of a picture shot just to prove 
      a person was in his hands -- think at the very poor quality
      pictures from real kidnappers and terrorists).
      \item
      As an expert photographer, he had to think that 
      the asymmetries in the picture would save him. 
      In particular
      \begin{enumerate}
      \item
      The picture shows an asymmetric disposition of the 
      furniture. Obviously he cannot tell which one is the correct
      one, but he could simply say that he was so sure it was 
      2:00 PM that, for example, the dresser had to be right of 
      fireplace and not on its left. He could simply require to check it.
      \item
      Finally, his wife wore a white rosette on her left. 
      This detail would allow him to claim with certainty 
      that the picture has been reversed (he knew how his wife
      was dressed, something that could be easily verified by the police,
      and, moreover, rosettes hang regularly left). 
      \end{enumerate}
      \end{enumerate}
      \label{fn:nota}
      }
      Then he rapidly goes to the 
      shelf, displaces one camera and, with no hesitation
      and no sign of doubt, 
      he picks up the one he used,
      that was visible, 
      but cleverly placed 
      in the back of others. Then, he opens that kind 
      of  old Polaroid-like camera and 
      shows the negative inside it 
      as the prove the picture was reversed
      and his alibi still valid.
\end{itemize}
But according to Columbo and his three colleagues, 
as well as to any TV-watcher, the {\it full action}
incriminates him. [Finally, although the confession
is irrelevant here,  he realizes his mistake
and his loss -- murderers in the Columbo series have 
usually some dignity.]

\subsection{How likely would have an innocent
person behave that way?}
As it easy to grasp, 
it is not just a question of picking up 
a camera out of 13. It is the entire sequence that 
is incompatible with an innocent person. 
Nobody, asked {\it directly} to pick up the camera 
used in a crime, would have done it on purpose, 
as a clear evidence of guilt. 

Certainly there was an 
{\it indirect} request, implicit in the Columbo
stratagem: ``find the negative''. 
But even 
an expert photographer would have not reacted that way
if he had been innocent. 

Let us assume it is reasonable he could overlook, 
in that particular,
dramatic moment, he had other ways out (see footnote  \ref{fn:nota})
and  only thought at the negative of the destroyed picture.
In this case he could have asked the 
policemen to take the camera 
and to look inside it.
Or he would have indicated the cameras behind Columbo's 
shoulder, suggesting that the negative {\it could}  
be in one of those cameras.\footnote{Nobody 
mentioned the camera was
in those shelfs or even in that room! 
(And TV watchers didn't get the 
information that 
Galesco knew that the camera was found by the police -- but
this could just be a minor detail.)
Moreover, only the killer and few policemen
knew that the negative was left inside it by the
murderer, a particular that is no obvious at all. As it was 
very improbable the killer used such an old-fashioned of 
camera.
Note in fact that the camera
was considered  a quite old one
already at the time the episode was set
and it was bought  in a second hand shop.
In fact I remember being wondering 
about that writer's choice,
until the very end:
it was done on the purpose, so that nobody but the
killer could think it was used to snap Mrs Galesco.
Clever!} 

An innocent person, even
put under dead stress and thinking that only the negative
of the destroyed picture
could save him, would perhaps jump
towards the shelf, take the first camera or the 
first few cameras he could reach and even
desperately shout
``look inside them!''. But he could have never
resolutely displaced other cameras, taken the correct
one on the back row and opened it, 
sure of finding the negative inside it. 

But not even a cool murderer would have reacted 
that way, 
as Galesco realized a bit too late.
The clever trick of Columbo was not only 
to ask indirectly the killer
to grasp the camera he used and that only he could  
recognize,
but, before that, to put him 
under stress in order to make him loose 
self control. 

%\newpage
\subsection{The verdict}
In summary, these are 
our reasonable beliefs that a person would have behaved 
that way producing that sequence of actions (${\cal A}$),
depending on whether he was a killer ($K$) or not ($\overline K$), 
maintaining or not self-control ($SC$/$\overline{SC}$):
\begin{description}
\item[$P({\cal A}\,|\,\overline K, I) = 0$]: 
this is the main point, that makes the hypothesis
innocent  definitely impossible.
\item[$P({\cal A}\,|\,K\cap {SC}, I) = 0$]: 
a cool murderer would have never reacted that way.
\item[$P({\cal A}\,|\,K\cap {\overline{SC}},I) > 0$]: 
this is the only hypotheses that can explain the action. 
Here `$>0$' stands for `not impossible', although not necessarily
`very probable'\footnote{Note 
that it is not required that one of the 
hypotheses should give with probability one, 
as it occurred instead of the toy example of section 
\ref{sec:1in13}. (See also Appendix G.)}. 
Let us say that, if Columbo
had planned his stratagem, based on a bluff, 
he knew there were some chances Galesco could reacted
that way, but he could not be sure about 
it.
\end{description}
Given this scenario, the `probabilistic inversion' 
is rather easy, as only one hypothesis remains possible:
that of a killer, who even had lost self-control.
%http://alligatographe.blogspot.com/2009/11/columbo-negative-reaction.html

%\newpage
\section{Comments and conclusions}
Well, this was meant to be a short note. Obviously
it is not just a comment to the New Scientist article,
that could have been contained in a couple of sentences.
In fact, discussing with several people, I felt that
yet another introduction to Bayesian reasoning,
not focused on physics issues, might be 
useful. 
So, at the end of the work, Columbo's cameras
were just an excuse. 

Let us now summarize what we have learned
and make further comments on some important issues.

First, we have to be aware that often we do not see 
`a fact' (e.g. Galesco killing his wife), 
but we 
infer it from other facts,
assuming a causal connection among them.\footnote{A quote 
by David Hume is in order (the subdivision
in paragraphs is mine):
\begin{quote}
{\sl
 All reasonings concerning matter of fact seem to be founded 
on the relation of Cause and Effect. 
By means of that relation alone we can go beyond the evidence 
of our memory and senses. 

If you were to ask a man, why he believes any matter of fact, 
which is absent; for instance, that his friend is in the country, 
or in France; he would give you a reason; 
and this reason would be some other fact; 
as a letter received from him, or the knowledge 
of his former resolutions and promises. 

A man finding a watch or any other machine in a desert island, 
would conclude that there had once been men in that island. 
All our reasonings concerning fact are of the same nature. 
And here it is constantly supposed that there is a connexion 
between the present fact and that which is inferred from it. 
Were there nothing to bind them together, the inference 
would be entirely precarious. 

The hearing of an articulate voice and rational discourse 
in the dark assures us of the presence of some person: 
Why? because these are the effects of the human make 
and fabric, and closely connected with it. 

If we anatomize all the other reasonings of this nature, 
we shall find that they are founded on the relation 
of cause and effect, and that this relation is either 
near or remote, direct or collateral.''
}\cite{ECHU}
\end{quote}
I would like to observe that too often we tend to take for
granted `a fact', forgetting that we didn't really observed
it, but we are {\it relying on a chain of testimonies and assumptions
that lead to it}. But some 
of them might fail (see footnote \ref{fn:Cromwell} and Appendix I). 
}
But sometimes the observed effect can be attributed to several
causes and, therefore, having observed an effect we
cannot be sure about its cause. Fortunately, since
our beliefs that each possible cause could produce
that effect are not equal, the observation 
modifies our beliefs on the different causes.
That is the essence of Bayesian reasoning. 
Since `Bayesian' has several 
flavors\footnote{Already in 1950 I.J. Good listed in Ref.~\cite{Good}
9 `theories of probability', some of which could be
called `Bayesian' and among which de Finetti's approach, 
just to make an example, does not appear.} 
in the literature,
I summarize the points of view expressed here:
\begin{itemize}
\item Probability simply states, in a quantitative way,
      how much we believe something. (If you like, 
      you can reason the other way around, thinking
      that something is highly improbable if you
      would be highly surprised if it occurs.\footnote{It is
      very interesting to observe how people are differently 
      surprised, in the sense of their emotional reaction, 
      depending on the occurrence of events that they considered
      more or less probable. Therefore, contrary to I.J. Good
      -- I have been a quite surprised about this -- according to 
      whom 
      ``to say that one degree of belief is more intense than another
      one is not intended to mean that there is more emotion attached
      to it''\,\cite{Good}, I am definitively closer to the position
      of Hume: 
      \begin{quote}
      {\sl 
       Nothing is more free than the imagination of man; and though it
cannot exceed that original stock of ideas furnished by the internal and
external senses, it has unlimited power of mixing, compounding,
separating, and dividing these ideas, in all the varieties of fiction
and vision. It can feign a train of events, with all the appearance of
reality, ascribe to them a particular time and place, conceive them as
existent, and paint them out to itself with every circumstance, that
belongs to any historical fact, which it believes with the greatest
certainty. Wherein, therefore, consists the difference between such a
fiction and belief? It lies not merely in any peculiar idea, which is
annexed to such a conception as commands our assent, and which is
wanting to every known fiction. For as the mind has authority over all
its ideas, it could voluntarily annex this particular idea to any
fiction, and consequently be able to believe whatever it pleases;
contrary to what we find by daily experience. We can, in our conception,
join the head of a man to the body of a horse; but it is not in our
power to believe that such an animal has ever really existed.

It follows, therefore, that the difference between {\it fiction} and
{\it belief} lies in some sentiment or feeling, which is annexed to the
latter, not to the former.
      }\cite{ECHU}
      \end{quote}
\label{fn:hume_prob}
      })
\item  ``Since the knowledge may be different with different persons or with
      the same person at different times, they may anticipate the same event
      with more or less confidence, and thus different numerical probabilities
      may be attached to the same event.''~\cite{SchroedingerA}
      This is the subjective nature of probability. 
\item Initial probabilities can be elicited, with all the vagueness
      of the case,\footnote{To state it in an explicit way,
      I admit, contrary to others, that probability values
      can be themselves uncertain, as discussed in footnote 
      \ref{note:mind_oscillations}. I understand that probabilistic
      statements about probability values might seem strange concepts
      (and this is the reason why I tried to avoid them in footnote
      \ref{note:mind_oscillations}), but I see nothing unnatural in 
      statements of the kind ``I am 50\% confidence
      that the expert will provide a value of probability in the
      range between 0.4 and 0.6'', as I would be ready to place 
      a 1:1 bet on the event that the quoted probability value will be 
      in that interval or outside it.} 
      on a pure subjective base (see Appendix C). 
      {\it Virtual bets} or comparisons with
      reference events can be useful `tools' to 
      force ourselves or experts to provide quantitative 
      statements of our/their
      beliefs. (See also Appendix C.)
\item Probabilities can (but need not) be evaluated by past frequencies and 
      can even be expressed in terms of expected 
      frequencies of `successes' in hypothetical trials.
      (See Appendix B.)
\item Probabilities of causes are not generated, but only 
      updated by new pieces of evidence.
\item Evidence is not only the `bare fact', but also all available
      information about it (see Appendix D). This point is 
      often overlooked, as in the criticisms to Columbo's episode
      raised by New Scientist~\cite{NS}. 
\item The update depends on how differently we believe that
      the various
      causes might produce the same effect (see also Appendix G).
\item The probability of a single hypothesis cannot be
      updated, if there isn't at least a second hypothesis
      to compare with, unless the hypothesis is absolutely
      incompatible with the effect [$P(E\,|\,H,I)=0$, and not
      `as little', for example, $10^{-9}$ or $10^{-23}$].
      Only in this special case an hypothesis is
      definitely falsified. (See Appendix G.)
\item In particular, if there is only one hypothesis 
      in the game, the final probability of this hypothesis
      will be one, no matter if it could produce the effect
      with very small probability (but not zero).
      \footnote{I have just learned from Ref.~\cite{Good} of the following
      Sherlock Holmes' principle: ``{\it If a hypothesis
      is initially very improbable but is the only one that 
      explains the facts,then it must be accepted}''. However,
      a few lines after, Good warns us that ``if the 
      only hypothesis that seems to explains the facts has very 
      small initial odds, then this is itself evidence that 
      some alternative hypotheses has been overlooked''\ldots}
\item Initial probabilities depend on the information 
      stored somehow in our brain; being, fortunately,
      each brain different from all others, it
      is quite natural to admit that, in lack of 
      `experimental data',``{\it quot capita, tot sententiae}''.
      (See Appendix C.)
\item In the probabilistic inference (i.e. that stems from 
      probability theory) the updating rule is univocally
      defined by Bayes' theorem (hence the adjective `Bayesian'
      related to these methods).
\item This objective updating rule makes final beliefs virtually
      independent from the initial ones, if rational people
      all share the same `solid' experimental information and 
      are ready to change their opinion (the latter
      disposition has been named {\it Cromwell's rule}
      by Dennis Lindley~\cite{Lindley}).
\item In the simple case that two hypotheses are involved,
      the most convenient way to express the Bayes' rule is 
      $$\mbox{final odds} = \mbox{Bayes factor} \times \mbox{initial odds},$$
      where the Bayes factor can be seen as the odds due to a 
      single piece of evidence, if the two hypotheses were considered 
      otherwise equally likely. (See also examples in Appendices F and G,
      as well as Appendix H, for comments 
      on statistical methods based on likelihood.)
\item In some cases -- almost always in scientific applications --
      Bayes factors can be calculated exactly, or almost exactly,
      in the sense that all experts will agree. 
      In many other real life cases their interpretation
      as `virtual' odds (in the sense stated above) allows to
      elicit them with the bet mechanism as any 
      subjective probability. (See Appendix C.)
\item Bayes factors due to several independent pieces of evidence
      factorize. 
\item The multiplicative updating rule can be turned into 
      an additive one using logarithms of the factors.
      (See Appendix E.)
\item The base 10 logarithms has been preferred here because
      they are easily related to the orders of magnitudes of the odds
      and the name  `judgement leanings' (JL) has been chosen to have no
      conflict with other terms already engaged in probability and statistics.
\item Each logarithmic addend has the meaning of weight of evidence,
      if the initial odds are taken as 0-th evidence. 
\item Individual contribution to the judgement 
      might be small in module and even 
      somehow uncertain, but, nevertheless, their combination
      might result into strong convincingness. (See Appendix G.)
\item In most real life cases there are not just
      two alternative causes and two possible effects. 
      Moreover, causes can be effects of other causes and effects 
      can be themselves causes of other effects. 
      All hypotheses in the game make up a complex `belief network'.
      Experts can certainly provide kinds of educated guesses
      to state how likely a cause can generate several effects,
      but the analysis of the full network goes well 
      beyond human capabilities, as discussed more extensively in 
      Appendix C and J. 
\item A next to simply case is when the evidence is mediated
      by a testimony. The formal treatment in Appendix I
      shows that, although experts can easily assess
      the required ingredients, the conclusions are really
      not so obvious. 
\item The question of the critical value of the 
      judgement leaning, above which a suspected
      can be condemned, goes beyond the purpose of this notes,
      focused on belief. That is a delicate decision problem that
      inherits all issues of assessing beliefs, to which 
      the evaluations of benefits and losses need to be added.
\end{itemize}
\vspace{0.5cm}\noindent
And Galesco? Come on, there is little to argue. \\
$[$Nevertheless,
the reading of the instructive New Scientist article is
warmly recommended!$]$

\vspace{1.0cm}\noindent
It is a pleasure to thank Pia and Maddalena, who
introduced me Columbo, and Dino Esposito, 
Paolo Agnoli
and Stefania Scaglia for having taken part to 
the post dinner jury that absolved him. 
The text has benefitted
of the careful reading by Dino, Paolo and Enrico Franco
(see in particular his interesting remark in footnote \ref{fn:enrico}).

%\mbox{}\\ \mbox{}\vspace{0.5cm}

%\mbox{}\\ \mbox{}\vspace{1.0cm}

\newpage
\appendix
\section{The rules of probability}
Let us summarize here the rules that degrees of belief
have to satisfy.
\vspace{-0.3cm}
\subsection{Basic rules}
Given any hypothesis $H$ (or $H_i$ if we have many of them),
also concerning the occurrence
of an event, and a given state of information $I$, probability 
assessments have to satisfy the following relations:
\begin{enumerate}
\item \fbox{$0\le P(H\,|\,I) \le 1$}
\item \fbox{$P(H\cup \overline H\,|\,I) = 1$}
\item \fbox{$P(H_i\cup H_j\,|\,I) = P(H_i\,|\,I) + P(H_j\,|\,I)$ \ \ if 
$H_i$ and $H_j$ cannot be true together}
\item  \fbox{$P(H_i\cap H_j\,|\,I) = P(H_i\,|\,H_j,I)\cdot P(H_j\,|\,I) 
= P(H_j\,|\,H_i,I)\cdot P(H_i\,|\,I)$}
\end{enumerate}
The first basic rule represents basically a conventional scale of 
probability, also indicated between 0 and 100\%. 

Basic rule 2 states that probability 1 is assigned to a 
{\it logical truth}, because either 
is true $H$ or its opposite (``{\it tertium non datur}''). 
Indeed $H\cup \overline H$
represent a logical, tautological certainty 
(a {\it tautology}, usually indicated with $\Omega$),
while $H\cap \overline H$ is a {\it contradiction},
that is something impossible, 
indicated by $\emptyset$.

The first three basic rules are also known the `axioms'
of probability,\footnote{Sometimes one hears of {\it axiomatic approach}
(or even {\it axiomatic interpretation} -- an expression
that in my opinion has very little sense)
of probability, also known as {\it axiomatic Kolmogorov approach}. 
In this approach `probabilities' are just real `numbers' in the range
$[0,1]$ that satisfy the axioms, 
with no interest on their {\it meaning}, 
i.e. {\it how they are perceived by the human mind}. 
This kind of approach
might be perfect for a pure mathematician, only interested to 
develop all mathematical consequences of the axioms. 
However it is not suited for applications, because, before we can use
the `numbers' resulting from such a probability theory, 
we have to understand what they mean. For this reason one 
might also hear
that ``probabilities are real numbers which obey the axioms
and that we need to `interpret' them'', an expression I deeply dislike.
I like much more the other way around: {\it probability is probability}
(how much we believe something) and probability
values {\it can be proved to obey} the four basic rules 
listed above, which can {\it then} considered by a pure mathematician 
the `axioms' from which a theory of probability can be built. 
}
while the inverses of the fourth one, e.g.
$ P(H_i\,|\,H_j,I)=P(H_i\cap H_j\,|\,I)/P(H_j\,|\,I)$, 
are called in most literature ``definition of conditional
probability''. In the approach followed here such a statement 
has no sense, because probability is always conditional probability
(note the ubiquitous `$I$' in all our formulae --
for further comments see section 10.3 of Ref.~\cite{BR}).
Note that when the condition $H_i$
does not change the probability of $H_j$, i.e.  
$P(H_i\,|\,H_j,I)= P(H_i\,|\,I)$, 
then  $H_i$ and $H_j$ are said to
be {\it independent in probability}. In this case the 
{\it joint probability} $P(H_i\cap H_j\,|\,I)$ is given 
by the so-called {\it product rule}, i.e. 
$P(H_i\cap H_j\,|\,I)=P(H_i\,|\,I)\cdot P(H_j\,|\,I)$.

These rules are automatically satisfied if probabilities
are evaluated from favorable over possible, equally probably cases. 
Also relative frequencies of occurrences in the past respect these
rules, with the little difference that the probabilistic 
interpretation of past relative frequencies is not really 
straightforward, as briefly discussed in the following appendix. 
That beliefs satisfy, in general,  the same basic rules can be 
proved in several ways. If we calibrate our degrees of beliefs 
against `standards', as illustrated in section 
\ref{sec:real_life}, this is quite easy to understand. 
Otherwise it can be proved by the normative principle of the 
{\it coherent bet}~\cite{BdF}.
%%\\\break\\ \mbox{}
\vspace{-0.3cm}
\subsection{Other important rules}
Important relations that follow from the basic rules are 
($A$ is also a generic hypothesis):
\begin{eqnarray}
P(\overline H\,|\,I) &=& 1 - P(H\,|\,I) \\
P(H\cap \overline H\,|\,I) &=& 0 \\
P(H_i\cup H_j\,|\,I) &=& P(H_i\,|\,I) + P(H_j\,|\,I) - P(H_i\cap H_j\,|\,I)
\label{eq:rul2}\\
P(A\,|\,I) &=& P(A\cap H\,|\,I) +  P(A\cap \overline H\,|\,I) \label{eq:rul3}\\
           &=& P(A\,|\,H,I)\cdot P(H\,|\,I) + P(A\,|\,\overline H,I)
\cdot P(\overline H\,|\,I) \label{eq:rul4} \\
P(A\,|\,I) &=& \sum_i P(A\cap H_i\,|\,I) 
\hspace{0.5cm}(\mbox{if $H_i$ form a {\it complete class}})
\label{eq:rul5}  \\
           &=& \sum_i P(A\,|\,H_i,I)\cdot P(H_i\,|\,I)
\hspace{0.5cm}(\mbox{idem})\,.  \label{eq:rul6}
\end{eqnarray}
The first two rules are quite obvious. Eq.~(\ref{eq:rul2})
is an extension of the third basic
rule in the case two hypotheses are not mutually exclusive.
In fact, if this is not case, the probability of 
$H_i\cap H_j$ is double counted and needs to be subtracted.
Eq.~(\ref{eq:rul3}) is also
very intuitive, because either $A$ is true together with $H$ or with
its opposite. 

Formally, Eq.~(\ref{eq:rul4}) follows from Eq.~(\ref{eq:rul3}) 
and basic rule 4. Its interpretation is that the probability of 
any hypothesis can be seen as `weighted average' of conditional
probabilities, with weights given by the probabilities of the 
conditionands [remember that $P(H\,|\,I)+P(\overline H\,|\,I)=1$
and therefore  Eq.~(\ref{eq:rul4})
can be rewritten as
$$
  P(A\,|\,I)= \frac{P(A\,|\,H,I)\cdot P(H\,|\,I) + P(A\,|\,\overline H,I)
\cdot P(\overline H\,|\,I)}{P(H\,|\,I)+P(\overline H\,|\,I)}\,,   
$$
that makes self evident its 
weighted average interpretation].

Eq.~(\ref{eq:rul5}) and (\ref{eq:rul6}) are simple extensions
of Eq.~(\ref{eq:rul3}) and (\ref{eq:rul4}) to a generic `complete class',
defined as a set of mutually exclusive hypotheses 
[$H_i\cap H_j=\emptyset$, i.e. $P(H_i\cap H_j\,|\,I)=0$], 
of which at least one must be true [$\cup_iH_i=\Omega$, 
i.e. $\sum_i P(H_i\,|\,I)=1$]. It follows then that Eq.~(\ref{eq:rul6})
can be rewritten as the (`more explicit') weighted average
$$
  P(A\,|\,I)=  \frac{\sum_i P(A\,|\,H_i,I)\cdot P(H_i\,|\,I)}
                         {\sum_i P(H_i\,|\,I)}\,.
$$
[Note that any hypothesis  $H$
and its opposite $\overline H$ form a complete class, 
because $P(H\cap \overline H\,|\,I)=0$ and 
$P(H\cup \overline H\,|\,I)=1$.]
%\\\break
%\newpage
\section{Belief versus frequency}
\subsection{Beliefs from past frequencies}
There is no doubt that 
\begin{quote}
{\sl
``where different effects
have been found to follow from causes, which are to {\it appearance}
exactly similar, all these various 
effects must occur to the mind in transferring the past to the future,
 and enter into our
consideration, when we determine the probability of the event. 
Though we give the preference to
that which has been found most usual, and believe that this 
effect will exist, we must not
overlook the other effects, but must assign to
each of them a particular weight and authority, in
proportion as we have found it to be more or less frequent.''  
}\cite{ECHU} 
\end{quote}
However, some comments about how our minds perform these 
operations are in order.\\
%\begin{enumerate}
%\item
\begin{minipage}{0.75\linewidth}
\vspace{.8mm}
\hspace{5.0mm}Before they are turned
into beliefs, observed frequencies are somehow smoothed, either 
intuitively or by mathematical algorithms. In both cases,
consciously or unconsciously, some models
of regularities are somehow `assumed' (a word that
in this context means exactly `believed').
Think, for example, at an experiment in which 
the number of counts are recorded in a defined interval of time,
under conditions apparently identical.  
Imagine that the numbers of counts in 20 independent measurements
are: 0, 0, 1, 0, 0, 0, 1, 2, 0, 0, 1, 1, 0, 4, 2, 0, 0, 0, 0, 1. 
The results are reported in the histogram. The question 
is ``what do we expect in the 21-st observation, provided the
experimental conditions remain unchanged?''. 
It is rather out of discussion that,
if a prize is offered on the occurrence of a count,
everyone will bet on 0, because 
it happened most
frequently. But
\vspace{1.1mm}
\end{minipage}
\begin{minipage}{0.25\linewidth}
\hspace{0.1\linewidth}\epsfig{file=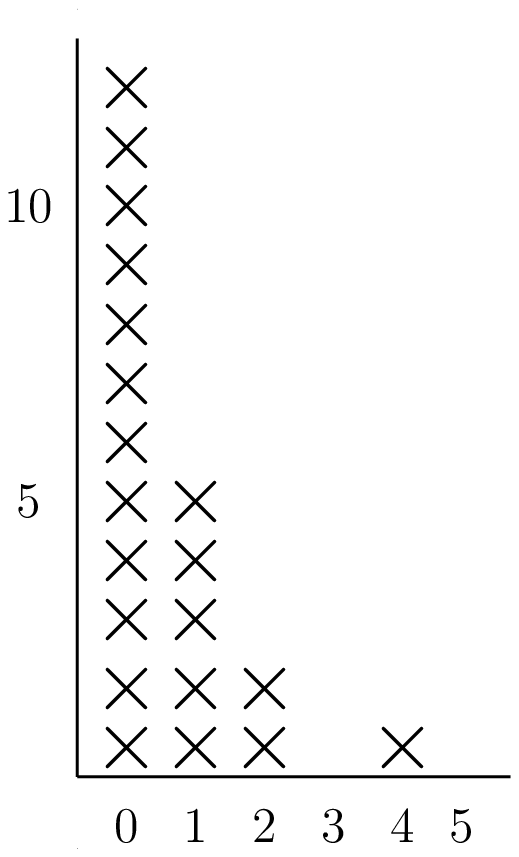,clip=,width=0.85\linewidth}
\end{minipage}
%\mbox{}\\
%\vspace{1.0mm}
 can we state that our belief is
 exactly
60\% (12/20)? Moreover, I am also pretty sure that, if you were
asked
 to place your bet on 3 or 4, you would prefer 3,
although this number of count has not occurred 
in the first 20 observations.
In an analogous way
 you might
 not believe that 5 is impossible. 
That is because we tend 
to see regularities in nature.\footnote{I find that
the following old joke conveys well the message. 
A philosopher, a physicist and a mathematician travel by train through
Scotland.
The train is going slowly and they see a cow walking along a country road 
parallel to the railway. The philosopher look at the others,
then very seriously states ``In Scotland cows are black''. 
The physicist replies that we cannot make such a generalization
from a single individual. We are only authorized to state, he maintains,
that ``In Scotland there is at least one black cow''. 
The mathematician looks well at cow,  
thinks a while, and then, he said, 
``I am afraid you are both incorrect. The most
we can say is that in Scotland at least one cow has a black side''.
}
Therefore going from past frequencies 
to probabilities can be quite a sophisticated process, that requires
a lot of assumptions (again {\it priors}!).

\subsection{Relative frequencies from beliefs}
The question of how relative frequencies of occurrence
follow from beliefs is much easier. It is
a simple consequence of probability theory and can 
be easily understood by anyone familiar with the binomial
distribution, taught in any elementary course on probability.
If we think at $n$ independent trials, for each of which 
we believe that the `success' will occur with probability $p$,
the {\it expected number} of successes is $np$, with 
a {\it standard uncertainty} $\sqrt{np(1-p)}$. 
We expect then a relative frequency $p$ [that is $(np)/n$]
with an uncertainty  $\sqrt{p(1-p)/n}$\, 
[that is $\sqrt{np(1-p)}/n$]. When $n$ is very large,
the uncertainty goes to zero and 
we become `practically sure' to observe a relative frequency
very close to $p$. This asymptotic feature goes under the name
of {\it Bernoulli theorem}. It is important to remark that 
this reasoning
can be {\it purely hypothetical} and has  nothing to do
with the so called frequentistic definition of probability.

To conclude this section, probabilities can be evaluated from 
(past) frequencies and (future, or hypothetical) 
frequencies can be evaluated from probabilities, but
{\it probability \underline{is not} frequency}.\footnote{The
following de Finetti's quote is in order.
%\begin{quote} 
{\sl
``For those who seek to connect the notion of probability 
with that of frequency,
results which relate probability  and frequency 
in some way (and especially those
results like  the `law of large numbers') play a pivotal r\^{o}le, 
providing support for the approach and for the identification
of the concepts. Logically speaking, however, one cannot escape 
from the dilemma posed by the fact that the same thing cannot both 
be assumed first as a definition and then proved as a theorem; nor
can one avoid the contradiction that arises from a definition which 
would  assume as certain something that the theorem only states 
to be very probable.''~\cite{BdF}
}
%\end{quote}
}
%\\\break\\ \mbox{}
%\\ \mbox{} \vspace{0.3cm}

%fine nuova A

%\subsection{Intuitions versus formal belief networks}\label{ss:intuition}
\section{Intuitions versus formal, possibly computer aided,
reasoning}
Contrary to `robotized Bayesians'\footnote{This expression
refers the robot of  E.T. Jaynes'~\cite{ETJ} and followers,
according to which probabilities should not be subjective.
Nevertheless, contrary to frequentists, they allow the possibility
of  `probability inversions' via Bayes' theorem, 
but they have difficulties with priors, that, according to them,
shouldn't be subjective. Their solution is that 
the evaluation of priors should be then
delegated to some `principles' 
(e.g. {\it Maximum Entropy} or {\it Jeffrey priors}). 
 But it is a matter of fact that 
unnecessary principles (that can be, anyway, used
as convenient rules in particular, well understood situations)
are easily misused (see e.g. comments on maximum
likelihood principle in the Appendix H 
-- several years ago, remarking this attitude by several
Bayesian fellows, I wrote a note on 
{\it Jeffreys priors versus experienced physicist priors; 
 arguments against objective Bayesian theory}, 
whose main contents went lately into Ref.~\cite{anxiety}), 
the approach
becomes dogmatic and uncritical use of some methods
might easily lead to absurd conclusions. 
For comments on 
anti-subjective criticisms 
(mainly those expressed in chapter 12
of Ref. \cite{ETJ}), see section 5 of Ref.~\cite{ME2000}.
As an example of a bizarre result, although considered 
by many Jaynes' followers as one of the jewels of their 
teacher's thought, let me mention the famous die problem.
``A die has been tossed a very large number $N$ of times, and we
are told that the average number of spots up per
toss was not 3.5, as we might expect from an honest die,
but 4.5. Translate this information into a probability 
assignment $P_n, n=1,2,\ldots,6$, for the $n$-th face
to come up on the next toss.''\cite{Jaynes} 
The celebrated Maximum Entropy solution is that the 
 probabilities for the six faces are,
in increasing order, 5.4\%, 7.9\%, 11.4\%, 18.5\%, 24.0\%
and  34.8\%. I have several times raised my perplexities
about the solution, but the reaction of Jaynes' followers
was, let's say, exaggerated.
Recently this result has been questioned 
by the somewhat quibbling Ref.~\cite{Mana} 
(one has to recognize that 
the original formulation of the problem had anyhow the assumption 
that the die was tossed a large number of times), 
which, however, also misses
the crucial point: {\it numbers on a die faces are just
labels}, having no intrinsic order, 
as instead it would be the case of 
the indications on a measuring device. 
I find absurd
making this kind of inferences without even giving a look 
at a real die! (Any reasonable person, used to try to observe
and understand nature,
 would first observe careful a die and try to guess 
how it could have been loaded to favor the faces having 
larger number of spots.)
\label{fn:Jaynes}}
I think it is quite natural that different
persons might have initially different opinions, 
that will necessarily influence the beliefs updated by 
experimental evidence, although the updating
rule is well defined, because based on probability
theory. But we have also seen, in a formal way,
that when the combined weight of evidence in favor of either
hypothesis is much larger than
the prior judgement leaning, i.e. 
$|\Delta\mbox{JL}_{1,2}(\mvec E,I)|\gg |\mbox{JL}_{1,2}(I)|$, 
then priors become irrelevant and we reach highly
inter-subjective conclusions. 

I am not in the position to try to discuss the internal processes
of the human mind that lead us to react in a certain way to
different {\it stimuli}. I only acknowledge that there are experts
of different fields that can make (in most case good) decisions in an 
fantastically short reacting time.
There is no need to think to doctors or engineer
in emergency situations, football players, fighter pilots,
and many other examples. It is enough to
observe us in the everyday actions
of driving a car or recognizing people from very partial 
information (and the context plays and important role! How many times
has happened to you not to immediately recognize/identify 
a neighbor, a waiter or a clerk
if you meet him/her in a place you didn't expect him/her at all?).
We are brought to think that much of the way in which external
information is processed is not analytical, but somehow 
hard-wired in the brain. 

A part of the automatic reasoning of the mind is innate, 
as we can understand observing children, animals, or even 
rational adults when they are possessed 
by pulsions and emotions. Another part comes from the experience
of the individual, where by `experience' it is meant 
all inputs received, of which he/she might be conscious
(like education and training) or unconscious, but
all processed and organized (again consciously or not) by the 
{\it causality principle}\cite{ECHU}, 
that allows us to anticipate (again consciously or not)
the consequences of our and somebody else's actions. 
%\footnote{Specie specific
%and individual specific: Darwin-Hume}
As a matter of fact, and coming to the main issue of this paper,
there is no doubt that
experienced policemen, lawyers and judges, thanks to their experience, 
have developed kinds of automatic reasonings, 
that we might call instinct or 
intuitive behavior (see footnote \ref{note:common-sense}) 
and that 
certainly help them in their work. 

We have seen in section \ref{sec:real_life}
that priors and even individual weights of evidence
can be elicited on a pure subjective way, possibly with 
the help of virtual bets or of comparison to reference events.
The problem arrives when the situation becomes a bit 
more complicate than just one cause and a couple of effects,
and the network of causes-effects become complex. 
Appendix I shows that the little complication of 
considering the possibility
that the evidence could be somehow reported in an 
erroneous way, as well known to psychologists,
of even fabricated by the investigators makes the problem
difficult and the intuition could fail. 
Appendix J shows an extension of the toy model of section 
\ref{sec:1in13}
in which several `testimonies' need to be taken into account.

In summary, the intuition of experts is fundamental to 
define the priors of the problem. It can be also very
important, and sometimes it is the only possibility
we have, to assess the degree
of belief that some causes can produce some effects, 
needed to evaluate the Bayes factors and, when the situation
becomes complex, to set up a `network of beliefs'
(see Appendix J). A different story 
is to process the resulting network, on the base
of the acquired evidences, in order to
evaluate the probabilities of interest. 
Intuition can be at lost, or miserably fail. 

\begin{figure}
\centering\epsfig{file=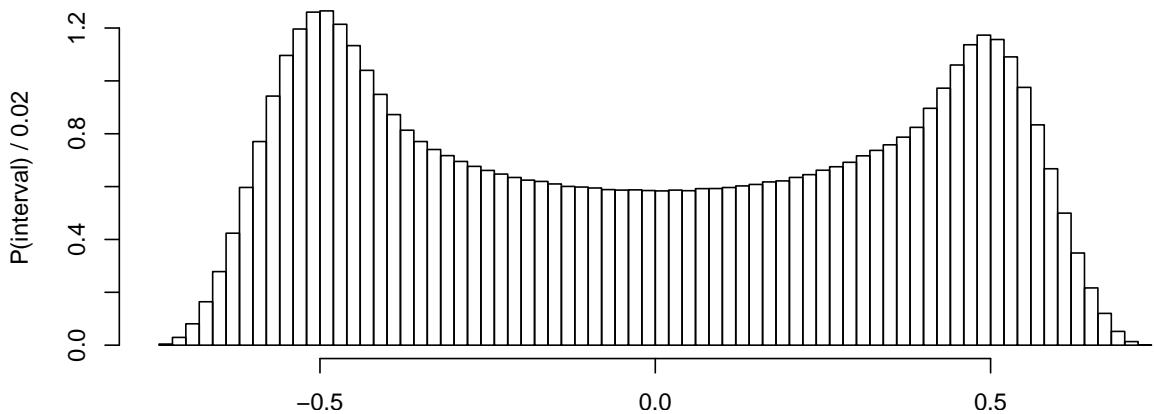,clip=,width=\linewidth}
\caption{\small \sf This histogram shows in a graphical 
way {\bf  `a' solution to the question at the end of Appendix C}
(details depend, obviously, on how the initial assumptions
have been modelled, but the gross features do not change 
if different reasonable models, consistent with the assumptions,
 are used -- here $x$ has been taken uniform between 9 and 11;
$y$ has been modelled with an asymmetric triangular distribution
ranging between 15 and 30, with maximum belief in 20).
\hspace{7.0cm}\mbox{ }
\mbox{ }\hspace{0.5cm}
The values of $z$ we \underline{have to} believe mostly are those 
around $-0.5$, but also all others in the range $-0.7$ to $0.7$
cannot be really neglected. In particularly, values around 0.5
are almost as likely as those around $-0.5$.
As we can see, there is about 50\% 
that $z$ occurs below 0 ($-0.02$, to be precise) 
and 50\% above. Note that, although the center of the 
distribution is around 0 ($-0.14$, to be precise), 
the most believable values are far from it. 
In other words, even if the {\it expected value} is $-0.14$
and the {\it standard uncertainty} (quantified by the 
{\it standard deviation}) is $0.40$, if a prize is assigned
to whom predicts the interval of width 0.02 in which 
the {\it uncertain number} $z$ will occur, 
we should place that interval at $-0.5$. 
%\hspace{13.0cm}\mbox{ }
%\mbox{ }\hspace{0.5cm}
Apart from the technical complications, the message of this
example is that one thing is to state the basic assumptions 
and subjective beliefs in some of the variables of 
the game, a much more complicate issue is to 
evaluate all {\it logical consequences of the premises}. 
In other words, if you agree on the premises of this problem, 
but not on the conclusions, you run into contradiction. 
Now, it is a matter of fact that {\it contradictions of this kind
are rather frequent} because the evaluation of the consequences
is not commonly done using formal logic and probability theory. 
The extension to complex belief networks is straightforward,
although, as we shall see in Appendix J, also a very simple
network is enough to challenge our ability to provide
intuitive answers.
}
\label{fig:dist_z}
\end{figure}

To make clearer the point consider this very rude example. 
Imagine you are interested in the variable $z$, that
you think for some reasons is
related to $x$ and $y$ by the following relation:
\begin{eqnarray*}
z &=& \frac{y\times\sin(\pi^4+x^2)}
           {\sqrt{x^3+y^2}}\,. 
\end{eqnarray*}
You might have good reason to state that $x$ 
is about 10, most likely not less than 9 and not more than 11,
and that in this interval you have no reasons to 
prefer a value with respect to another one. 
Similarly, you might thing that the value of $y$ you trust mostly
is 20, but it could go down to 15 and up to 30 with decreasing beliefs.
What do you expect for $z$.  
Which values of $z$ should you believe, consistently with your
basic assumptions? If a rich prize is give to the person
that predicts the interval of width 0.02 in which $z$ will occur, 
which interval would you choose?
What is the value of $z$ (let us 
call it $z_m$) such that there
is 50\% chance that $z$ will occur below this value?
What is the probability that $z$ will be above 10?
[The solution is in  next page (figure \ref{fig:dist_z}).]

Anyway, if you consider this example a bit too `technical'
you might want to check the capabilities of 
your intuition on the much 
simpler one of Appendix J. 
(Try first to read  the caption of figure 
\ref{fig:bn_mon_0} and to reply the questions.)
\\ \mbox{} \vspace{1.0cm}

%\newpage

%\mbox{}

%\mbox{}

%\newpage

\section{Bare facts and complete state of information}
As it has been extensively discussed in section \ref{sec:finale}, saying that 
a person has taken a camera out of thirteen is a piece of
information, but it is not all, and it is not enough to
update correctly our beliefs.

 This is true in general,
even in fields of research that are considered by outsider
to be the realm of objectivity, where only `facts' count.
Stated with Peter Galison words~\cite{Galison},
\begin{quote}
{\sl
``Experiments begin and end in a matrix of beliefs. 
\ldots beliefs in instrument type, in programs of experiment
enquiry, in the trained, individual judgments about every local behavior
of pieces of apparatus.''
}\\
$[$Then,  taking as an example the discovery of the positron:$]$\\
{\sl``Taken out of time there is no sense to
the judgment that Anderson's
 track 75
is a positive electron;
its textbook reproduction has been denuded of the prior experience
 that made Anderson
 confident in the cloud chamber, the magnet,
 the optics, and the photography.''
}
\end{quote}
My preferred toy examples to convey this important messages
are {\it the three box problem(s)'} and
{\it the two envelopes `paradox'} (see section 3.13 of Ref.~\cite{BR} --
I remind briefly here only the box ones).
The box problems are a series of recreational/educational 
problems, the basic
one being rather famous as `Monthy Hall problem'. 
The great majority of people (my usual target are physics
PhD students) get mad with them because they have not
been educated to take into account all available information.
Therefore they have quite some difficulties to understand 
that if a contestant has taken one box (yet unopened) 
and there is still another un-opened box to choose, the 
probability that this box contains the prize (only one of the three
boxes does) depends on whether 
the opened (and empty) box was got by chance or was chosen with 
the intention to take a box without 
prize.\footnote{Reading the draft of this paper, my colleague
Enrico Franco has remarked 
that in the way the box problems (or the Monthy Hall) are
presented there are additional pieces of information which are usually
neglected, as I also did in Ref.~\cite{BR}
(`then' was not underlined in the original):
\begin{quote}
(1) In the first case, imagine two contestants, each of whom
 chooses one box at random. Contestant $B$ opens his chosen 
box and finds it does not contain the prize. \underline{Then} the
presenter offers player $A$ the opportunity to exchange his
box, still un-opened, with the third box. \ldots
\\
(2) In the second case there is only one contestant, $A$. 
After he has chosen one box the presenter tells him that, 
although the boxes are identical, he knows which one contains the prize. 
\underline{Then} he says that, out of the two remaining boxes, 
he will open  one that does not
contain the prize.\ldots~\cite{BR}
\end{quote}
It makes quite some difference if the conductor announces
he will
propose the exchange
\underline{before} the boxe(s) is/are initially
taken by the contestant(s) that, or if he does it \underline{later}, as I usually
formulate the problems. In the latter case, in fact, contestant $A$ 
can have a legitimate doubt concerning the malicious intention of the 
conductor, who might want to induce him to lose. 
Mathematics oriented guys
would argue then that the problem does have a solution. 
But the question is that in real life 
one \underline{has to} act, 
and one has to finally make his decision, based on the best
knowledge of the game and of the conductor, 
in a finite amount of time.\label{fn:enrico}}
 (And there is often
somebody in the audience that when he/she listens the 
formulation of the problem in which the box was opened by chance,
he/she smiles at the others, and than gives the solution\ldots of the
version in which the conductor opens on purpose and empty box.)

I found that the issue of considering 
into account  all available information
is shown in a particular convincing way 
in the `three prisoner paradox' (isomorph\footnote{This is true
only neglecting the complication taken into account 
in the previous footnote. Indeed, in one case 
the `exchange game' is initiated by the conductor, 
while in the second by the prisoner, 
therefore Enrico Franco's comment does not apply to
the three prisoner problem.} to Monthy Hall,
but more a headache than this, perhaps because it involves humans)
and in the `thousand prisoner problem' of 
Ref.~\cite{Pearl}: 
not only bare facts enter the evaluation
of probability, but also all contextual knowledge about them,
including the question asked to acquire their knowledge.
%\break \\ \mbox{} \vspace{2cm} \mbox{}

%\newpage
\section{Some remarks on the use of logarithmic
updating of the odds}
The idea of using (natural) logarithms of the odds 
is quite old, going back, as far as I know,
to Charles Sanders Peirce~\cite{Peirce}.
He related them to what he called
{\it feeling of belief} (or {\it intensity of belief}), 
that, according to him,
 ``should be as the logarithm 
of the chance, this latter being the expression of the state
of facts which produces the belief''~\cite{Peirce}, where by `chance'
he meant exactly probability ratios, i.e. the odds.

Peirce proposed his
''thermometer for the proper intensity of belief''~\cite{Peirce} 
for several reasons. 
\begin{itemize}
\item
First because of 
considerations that when the odds
go to zero or to infinity, then the intensity of belief on either
hypothesis goes to infinity;\footnote{In  this respect,
 belief becomes similar
to other human sentiments, for which in normal speech we use a scale that
goes to infinity -- think at expressions like `infinite love', 
`infinite hate', and so on
(see also footnote \ref{fn:hume_prob}).
}  
when ``an even chance is reached
[the feeling of believing] should completely vanish
and not incline either toward or away from the proposition.''~\cite{Peirce}
The logarithmic function is the simplest one to achieve the
desired feature. (Another interesting feature of the odds is
described in footnote \ref{fn:logP}.) 
\item
Then because 
(expressing the question in our terms), 
if we started from a state of indifference (initial odds equal to 1),
each piece of evidence should produce odds equal
to its Bayes factor [our $\tilde O_{i,j}(E_i)$]. The combined odds will be  
the product of the individual odds 
[Eq.~\ref{eq:product_Odds}].
But, mixing now Pierce's and our terminology, 
when we combine several arguments (pieces of evidence),
they ``ought to produce a belief equal to the sum
of the intensities of belief which either would produce
separately''.~\cite{Peirce} Then ``because we have seen that the chances
of independent concurrent arguments are to be multiplied together
to get the chance of their combination, and therefore
the quantities which best express the intensities
of belief should be
such that they are to be {\it added} when the {\it chances}
are multiplied\ldots Now, the logarithm of the chance
is the only quantity which fulfills this condition''.~\cite{Peirce}
\item
Finally, Peirce justifies his choice by the fact
that human perceptions go often as the logarithm of
the stimulus (think at subjective feeling of 
sound and light -- even `utility', 
meant as the `value of money' is supposed 
to grow logarithmically
with the amount of money): 
``There is a general law of sensibility,
called Fechner's psychophysical law. It is that the intensity
of any sensation is proportional to the logarithm of the external 
force which produces it.''\cite{Peirce}
(Table \ref{tab:JL_BF_P} provides a comparisons between the different 
quantities involved, to show that the human sensitivity on probabilistic 
judgement is indeed logarithmic, with a resolution about the first
decimal digit of the base 10 logarithms.)
\end{itemize}

As far as the logarithms in question, I have
done a short research on their use,
which, actually, lead me to discover Peirce's
{\it Probability of Intuition}~\cite{Peirce}
and Good's {\it Probability and the weighing of 
Evidence}~\cite{Good}.\footnote{Peirce article is a mix of 
interesting intuitions and confused arguments, as in the 
``bag of beans'' example of pages 709-710
(he does not understand the difference
between the observation of 20 black beans 
and that of 1010 black and 990 white for the 
evaluation of the probability that another bean extracted
from the same bag is white or black, arriving thus
to a kind of paradox -- from Bayes' rule it is clear
that weights of evidence sum up to form the intensity
of belief on two bag compositions, not 
on the outcomes from the boxes~\cite{Peirce_beans}). Of a different class is 
Good's book, one of the best on probabilistic reasoning
I have met so far, perhaps because I feel myself
often in tune with {\it Good thinking} (including the passion
for footnotes and internal cross references shown in Ref.~\cite{Good}).}
As far as I understand, without pretension
of completeness or historical exactness:
\begin{itemize}
\item Peirce' `chances' are introduced as if they were our odds,
but are used if they were \underline{Bayes factors} 
(``the chances of 
independent concurrent arguments are to be multiplied together
to get the chance of their combination''~\cite{Peirce}). 
Then he takes the {\it natural} logarithm of these `chances',
to which he also associates an idea of {\it weight of evidence}
(``our belief ought to be proportional to the
weight of evidence, in the sense, that two arguments which are
entirely independent, neither weakening nor strengthening each other,
ought, when they concur, to produce a belief equal to the sum of the
intensities of belief which either would produce 
separately''~\cite{Peirce}).
\item According to Ref.~\cite{MT} 
the modern use of the logarithms of the
\underline{odds} seem to go back to I.J. Good, who used
to call {\it log-odds} the  {\it natural} logarithm
of the odds.\footnote{But Goods mentions that ``In 1936 Jeffreys
had already appreciated the importance of the logarithm
of the [Bayes] factor and had suggested for it the name 
`support'.''~\cite{Good}}
\item However, reading later Ref.~\cite{MT} it is clear that
      Good, following a suggestion of A.M. Turing, 
      proposes a decibel-like (db) 
      notation\footnote{``In acoustic and electrical engineering
      the bel is the logarithm to base 10 of the ratio of
      two intensities of sound. Similarly, if $f$ is the [Bayes] factor
      in favor of a hypothesis has gained 
      $\log_{10}f$ bels, or $(10\,\log_{10}f)$\,db.''~\cite{Good}
      [Good uses the name `factor' for what we call Bayes factor, 
       ``the factor by which the initial odds of $H$
       must be multiplied in order to obtain the final odds. 
       Dr. A.M. Turing suggested in a conversation in 1940
       that the word `factor' should be regarded as the technical
       term in this connexion, and that it could be more fully
       described as {\it the factor in favor of the hypothesis
       $H$ in virtue of the result of the experiment}.''~\cite{Good}]}, 
       giving proper names both to the logarithm of the 
       odds and to the logarithm of the Bayes factor:
       \begin{itemize}
       \item ``$(10\,\log_{10}f)$\,db \ldots 
             may be also described as the {\it weight of 
             evidence} or amount of information for $H$ 
             given $E$''~\cite{Good};
       \item ``$(10\,\log_{10}o)$\,db may be called the {\it plausibility}
             corresponding to odds $o$''~\cite{Good}. 
       \end{itemize}
       It follows then that 
       \begin{eqnarray}
       \mbox{``Plausibility gained} &=& \mbox{weight of evidence''}.
       \label{eq:plausibility_evidence}\cite{Good}
       \end{eqnarray} 
\item Decibel-like logarithms of the \underline{odds} are used 
      since at least forty years with under the name 
      {\it evidence}.~\cite{Jaynes}\,. 
\end{itemize}
Personally, I think that 
the decibel-like definition
is not very essential (decibels themselves
tend already to confuse normal people,
also because for some applications
the factor 10 is replaced by a factor 20).
Instead, as far as names are concerned:
\begin{itemize}
\item `plausibility' is difficult to defend, because 
      it is too similar to probability in everyday use,
      and, as far as I understand, has decayed;
\item `weight of evidence' seems to be a good choice,
       for the reasons already well clear to Peirce.
\item `evidence' in the sense of Ref.~\cite{Jaynes}
       seems, instead, quite bad for a couple of reasons:
\begin{itemize}
\item First, 
 because `evidence'
has already too many meanings, including, 
in the Bayesian literature, the denominator
of the r.h.s. of Eq.~(\ref{eq:PHi|E}).
\item
Second, because this name is given to the logs of the odds
(including the initial ones), but not to those of the 
Bayes factors to which no name is given. Therefore, 
the name `evidence', as used in Ref.~\cite{Jaynes}
in this context, is not related to the evidence.
\end{itemize}
\end{itemize}
I have taken the liberty to use the expression `judgment
leaning'
first because it evokes the famous {\it balance of Justice},
then because all other expressions I thought about have already 
a specific meaning, and some of them even several 
meanings.\footnote{Many controversies in 
probability and statistics arise because there
is no agreement on the meaning of the words
(including `probability' and `statistics'), 
or because some refuse to accept this 
fact. For example, I am perfectly aware
that many people, especially my friends physicists,
tend to to assign to the word `probability' the meaning
of a kind of {\it propension} `nature' has to behave 
more in a particular way than in other way, 
although in many other cases -- and more often!  -- 
they also mean by the same word
how much they believe
something (see e.g. chapters 1 and 10 of Ref.~\cite{BR}).
For example, one might like to think that kind $B_1$
boxes of section \ref{sec:1in13} \underline{have} a 100\% 
{\it propensity} to produce white balls and 0
to produce black balls, while type $B_2$
\underline{have} 7.7\% propension to produce white
and 92.3\% to produce black. Therefore, if one
knows the box composition and 
is only interested to the outcome of the extraction, then
probability and propensity coincide in value. But if the composition
is unknown this is no longer true, as we shall see
in Appendix J. [By the way, all interesting
questions we shall see in Appendix J have no meaning (and no
clean answers) for ideologizied guy who refuse to accept
that probability primarily means how much we believe something.
(See also comments in Appendix H.)]
\label{fn:propension}
}
It is clear, especially comparing 
Eq.~(\ref{eq:plausibility_evidence}) 
with Eq.~(\ref{eq:Delta_JL_diffJL}),
that, besides the factor ten multiplying
the base ten logarithms and the notation,
I am quite in tune with Good. 
I have also to admit I like Peirce'
`intensity of belief' to name the JL's,
although it is too similar to `degree of belief',
already widely used to mean something else.

So, in summary, these are the symbols and names used here:
\begin{description}
\item[$\mbox{JL}_{i,j}(\cdot)$] is the 
{\it judgement leaning} in favor
of hypothesis $i$ and against $j$, with the conditions in parenthesis.
If we only consider
an hypothesis ($H$) and its opposite $\overline H$, 
that could be possibly related to the occurrence of the
event $E$ or its opposite $\overline E$, also the notation
 $\mbox{JL}_{H}(\cdot)$, or  $\mbox{JL}_{E}(\cdot)$, will
be used (as   in table \ref{tab:JL_ET} of Appendix I).\\
(Sometimes I have also tempted to call a JL `intensity
of belief' if it is clear from the contest that the expression
does not refer to a probability.)
\item[$\Delta\mbox{JL}_{i,j}(\cdot)$], with the
same meaning of the subscript and of the argument, is the
variation of judgement leaning produced by 
a piece of evidence and it is called here 
{\it weight of evidence}, although it differs by a factor from 
the analogous names used by 
Peirce and Good\footnote{$\log_{10}x =  \ln x /\ln 10
= (10\,\log_{10}x)/10$.}.
\end{description}

%\\ \mbox{} \vspace{2.0cm}
%\newpage
\section{AIDS test}
Let us make an example 
of general interest, that exhibits some of the
issues that also arise in forensics.

Imagine an Italian citizen is chosen {\it at random} 
to undergo an AIDS test. Let us assume
the analysis used to test for HIV infection
is not perfect. In particular, infected people ($\mbox{HIV}$)
are declared `positive' ($\mbox{Pos}$) with 99.9\% probability
and `negative' ($\mbox{Neg}$) with 0.1\%; 
there is, instead, a 0.2\% chance
a healthy person 
($\overline{\mbox{HIV}}$) is told positive 
(and 99.8\% negative). 

The other information we need is the prevalence of the virus
in Italy, from which we evaluate our initial belief that 
the randomly chosen person is infect. We take  1/400 or 0.25\%
(roughly 150 thousands in a population of 60 millions). 

\begin{figure}
\centering\epsfig{file=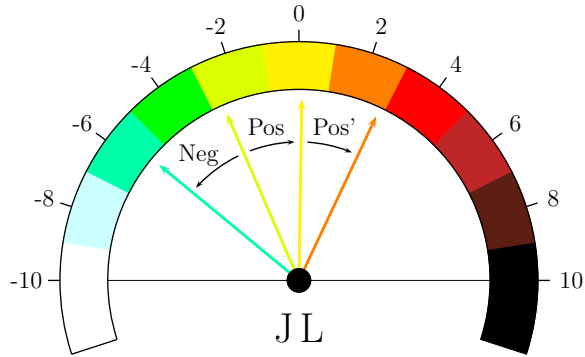,clip=,width=0.5\linewidth}
\caption{\small \sf AIDS test illustrated with judgement leanings.}
\label{fig:jl_aids}
\end{figure}

To summarize, these are the pieces of information relevant
to work the exercise:\footnote{The performance of the test
are of pure fantasy, while the prevalence is somehow
realistic, although not pretended to be the real one.
But it will be clear that the result is rather insensitive
on the precise figures.}
\begin{eqnarray*}
P(\mbox{Pos}\,|\,\mbox{HIV},I) & = & 99.9\%,  \\
P(\mbox{Neg}\,|\,\mbox{HIV},I) & = & 0.1\%,  \\
P(\mbox{Pos}\,|\,\overline{\mbox{HIV}},I) & = & 0.2\%\\
P(\mbox{Neg}\,|\,\overline{\mbox{HIV}},I) & = & 99.8\%\\
P(\mbox{HIV}\,|\,I) & = & 0.25\% \\
P(\overline{\mbox{HIV}}\,|\,I) & = & 99.75\%,
\end{eqnarray*}
from which we can calculate initial odds, Bayes factors
and JL's [we use here the notation $O_{\mbox{{\footnotesize HIV}}}(I)$, 
instead of our usual $O_{1,2}(I)$ to indicate odds in favor of 
the hypothesis HIV
and against the opposite hypothesis ($\overline{\mbox{HIV}}$); 
similarly for $\mbox{JL}_{HIV}$ and $\Delta\mbox{JL}_{HIV}$]:
$$
\begin{array}{rclcl}
O_{\mbox{{\footnotesize HIV}}}(I) &=& 1/399 = 0.0025 
\ \ \ &\Rightarrow& \mbox{JL}_{\mbox{{\footnotesize HIV}}}(I)= -2.6 \\
\tilde O_{\mbox{{\footnotesize HIV}}}(\mbox{Pos},I) &=& 99.9/0.2 = 499.5 
\ \ \  &\Rightarrow& \Delta\mbox{JL}_{\mbox{{\footnotesize HIV}}}(\mbox{Pos},I) = +2.7\\
\tilde O_{\mbox{{\footnotesize HIV}}}(\mbox{Neg},I) &=& 0.1/99.8 = 1/998 = 0.001002 
\ \ \  &\Rightarrow& \Delta\mbox{JL}_{\mbox{{\footnotesize HIV}}}(\mbox{Neg},I) = -3.0\,.
\end{array}
$$
A positive result adds a weight of evidence of 2.7 to $-2.6$, 
yielding the negligible leaning of $+0.1$. Instead a negative 
result has the negative weight of $-3.0$, shifting the leaning
to $-5.6$, definitely on the safe side (see fig. \ref{fig:jl_aids}).

The figure shows also the effect of a second, 
{\it independent}\footnote{Note that `independent' does not mean
the analysis has simply been done by somebody else,
possibly in a different laboratory, but also that the {\it principle 
of measurement} is independent.}
analysis, having the same performances of the first one and
in which the person results again positive. As it clear from
the figure, the same conclusion would be reached if only one
test was done on a subject for which a doctor could be in serious
doubt if he/she had AIDS or not ($\mbox{JL}\approx 0$).

From this little example we learn that if we want to have 
a good discrimination power of a test, it should have 
a $\Delta$JL very large in module. Absolute discrimination can only
be achieved if the weight of evidence is infinite, i.e. if
either hypothesis is impossible given the observation.
%\\ \mbox{} \vspace{3.0cm}
%\mbox{}

\section{Which generator?} 
\begin{figure}
\centering\epsfig{file=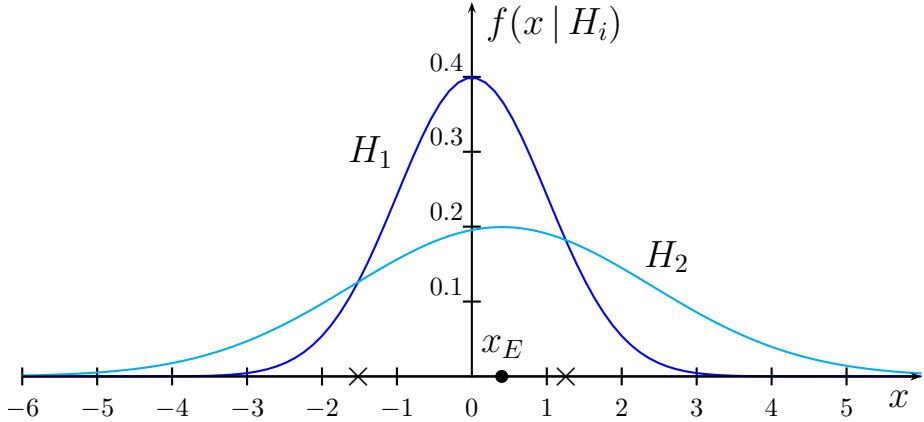,clip=,width=0.8\linewidth}
\caption{\small \sf Which random number generator has produced $x_E$?
Which hypothesis favors the points indicated by `$\times$'?}
\label{fig:two_normal}
\end{figure}

%(Combining pieces of evidence
%non uniform
% varying in size and even 
%sometimes contradictory pieces of information)}
%(The small probability of the observation is irrelevant!)}
Imagine two (pseudo-) random number generators: $H_1$,
Gaussian with mean 0 and standard deviation 1,
and $H_2$, also Gaussian, but
with mean 0.4 and standard deviation 2 (see figure \ref{fig:two_normal}).

A program chooses at random, with equal probability,
$H_1$ or $H_2$; then the generator produces a number,
that, rounded to the 7-th decimal digit, is
$x_E=0.3986964$. The question is, from which random generator
does $x_E$ come from? 

At this point, 
the problem is rather easy to solve, if we know the probability
of each generator to give $x_E$. They are\footnote{The curves
$f(x\,|\,H_i)$
 in figure \ref{fig:two_normal}  
 represent {\it probability density functions} (`pdf'), 
i.e. they give the probability per unit $x$, i.e. 
$P([x-\Delta x/2,x+\Delta x/2])/\Delta x$,
for small $\Delta x$ (remember that `densities' are always local). 
Rounding to the 7-th digit means that the number before 
rounding was in the interval of $\Delta x=10^{-7}$ 
centered $x_E$. It follows that the probability a generator
would produce that number can be calculated as 
$f(x_E\,|\,H_i)\times \Delta x$. Indeed, we can see
that in the calculation of Bayes factors the width $\Delta x$ simplifies
and what really matter is the ratio of the two pdf's, i.e.
\begin{eqnarray*}
\tilde O_{1,2}(x_E,I)&=&
\frac{P(x_E\,|\,H_1)}{P(x_E\,|\,H_2)} = 
\frac{f(x_E\,|\,H_1)\times\Delta x}{f(x_E\,|\,H_2)\times\Delta x}\,
= \frac{f(x_E\,|\,H_1)}{f(x_E\,|\,H_2)}\,. 
\end{eqnarray*}
The Bayes factor is therefore the ratio
of the ordinates of the curves in figure  \ref{fig:two_normal} 
for the same $x_E$. Note that 
$f(x_E\,|\,H_1)\times \Delta x$ can be small at will, 
but, nevertheless, hypothesis $H_1$ can receive  
a very high weight of evidence from $x_E$
if $f(x_E\,|\,H_1) \gg f(x_E\,|\,H_2)$.
}
\begin{eqnarray*}
P(x_E\,|\,H_1,I) &=& 3.68\times 10^{-8} \ \ \
(\mbox{1 in $\approx$\,27 millions}) \\
P(x_E\,|\,H_2,I) &=& 1.99\times 10^{-8} \ \ \ 
(\mbox{1 in $\approx$\,50 millions})\,,
\end{eqnarray*}
from which we can calculate Bayes factor and weight of evidence:
$$
\begin{array}{rclcl}
\tilde O_{1,2}(x_E,I) &=& 1.85
\ \ \  &\Rightarrow& \Delta\mbox{JL}_{1,2}(x_E,I) = +0.27\,.
\end{array}
$$
Therefore, the observation of $x_E$ provides a slight evidence
in favor of $H_1$, {\it no matter if 
this generator has very little probability to give $x_E$},
as it has very little probability to give {\it any} particular
number.

What matters when comparing hypotheses is never, 
stated in general terms,
the absolute probability $P(E\,|\,H_i,I)$.
In particular, it doesn't make sense saying
 ``$P(H_i\,|\,E,I)$ is small because $P(E\,|\,H_i,I)$ is 
small''.\footnote{Sometimes this might be
qualitatively correct,  because it easy to {\it imagine
there could be}
an alternative hypothesis $H_j$ such that: 
\begin{enumerate}
\item  $P(E\,|\,H_j,I) \gg P(E\,|\,H_i,I)$, such that 
       the Bayes factor is strongly in favor of $H_j$;
\item  $P(H_j\,|\,I)\approx P(H_i\,|\,I) $, that is $H_j$ 
       is roughly as credible as  $H_i$. 
\end{enumerate}
(For details see section 10.8 of Ref.\cite{BR}.)
} 
As a consequence, from a consistent probabilistic point
of view, {\it it makes no sense to test a 
single, isolated hypothesis}, using 
`funny arguments',
like how far if $x_E$ from the peak of $f(x\,|\,H_i)$,
or how large is the area below  $f(x\,|\,H_i)$ from $x=x_E$
to infinity. In particular, if two models 
give exactly the same probability to produce an observation,
like the two points indicated by `$\times$' in fig. 
\ref{fig:two_normal}, the evidence provided by
this observation is absolutely irrelevant 
[$\Delta\mbox{JL}_{1,2}(\mbox{`$\times$'})=0$].

To get a bit familiar with the weight of evidence in favor of either
hypothesis provided by different observations, the following table, 
reporting Bayes factors and JL's due to the integers between $-6$ and $+6$,
might be useful.\\
%\begin{center}
\begin{center}
\begin{tabular}{rclr}
\multicolumn{1}{c}{$x_E$} && \multicolumn{1}{c}{$\tilde O_{1,2}(x_E)$} & 
\multicolumn{1}{c}{\ \ $\Delta$JL$_{1,2}(x_E)$} \\
\hline 
$-6$ && $5.1\times 10^{-6}$ & $-5.3$ \\
$-5$ && $2.9\times 10^{-4}$ & $-3.5$ \\
$-4$ && $7.5\times 10^{-3}$ & $-2.1$ \\
$-3$ && $9.4\times 10^{-2}$ & $-1.0$ \\
$-2$ && $0.56$              & $-0.3$ \\
$-1$ && $1.5$               & $ 0.2$ \\
$0$  && $2.0$               & $ 0.3$ \\
$1$  && $1.3$               & $ 0.1$ \\
$2$  && $0.37$              & $-0.4$ \\ 
$3$  && $5.2\times 10^{-2}$ & $-1.3$ \\
$4$  && $3.4\times 10^{-3}$ & $-2.5$ \\
$5$  && $1.0\times 10^{-4}$ & $-4.0$ \\
$6$  && $1.5\times 10^{-6}$ & $-5.8$
\end{tabular}
%\caption{\small \sf Evidence produced by some values of $x$.}
%\label{fig:x_H1_H2}
\end{center}
%\end{table}
As we see from this table, and as we better understand from 
figure \ref{fig:two_normal}, numbers large in module are 
in favor of $H_2$, and very large ones are in its strong favor.
Instead, the numbers laying in the interval defined by the two 
points marked  in the figure by a cross provide evidence in favor of
$H_1$. However, while individual pieces of evidence in favor of $H_1$
can only be weak (the maximum of $\Delta$JL is about 0.3, reached
around $x=0$, namely $-0.13$, to be precise, 
for which $\Delta$JL reaches 0.313), 
those in favor of the alternative hypothesis can 
be sometimes very large. It follows then that one 
gets easier convinced of $H_2$ rather than of $H_1$.

We can check this by a little simulation. We choose
a model, extract 50 random variables and analyze
the data as if we didn't know which generator produced them,
although
considering $H_1$ and $H_2$ equally likely. We expect 
that, as we go on with the extractions, the pieces 
of evidence accumulate until we {\it possibly}
reach a level of practical certainty.
Obviously, the individual pieces
of evidence do not provide the same $\Delta$JL, and also the sign
can fluctuate, although we expect more positive contributions
if the points are generated by $H_1$ and the other way around
if they came from $H_2$. 
 Therefore, as a function 
of the number of extractions the accumulated weight 
of evidence follows a kind of {\it asymmetric random walk}
(imagine the JL indicator fluctuating 
as the simulated experiment goes on, but drifting
`in average'
 in one direction). 
\begin{figure}
\centering\epsfig{file=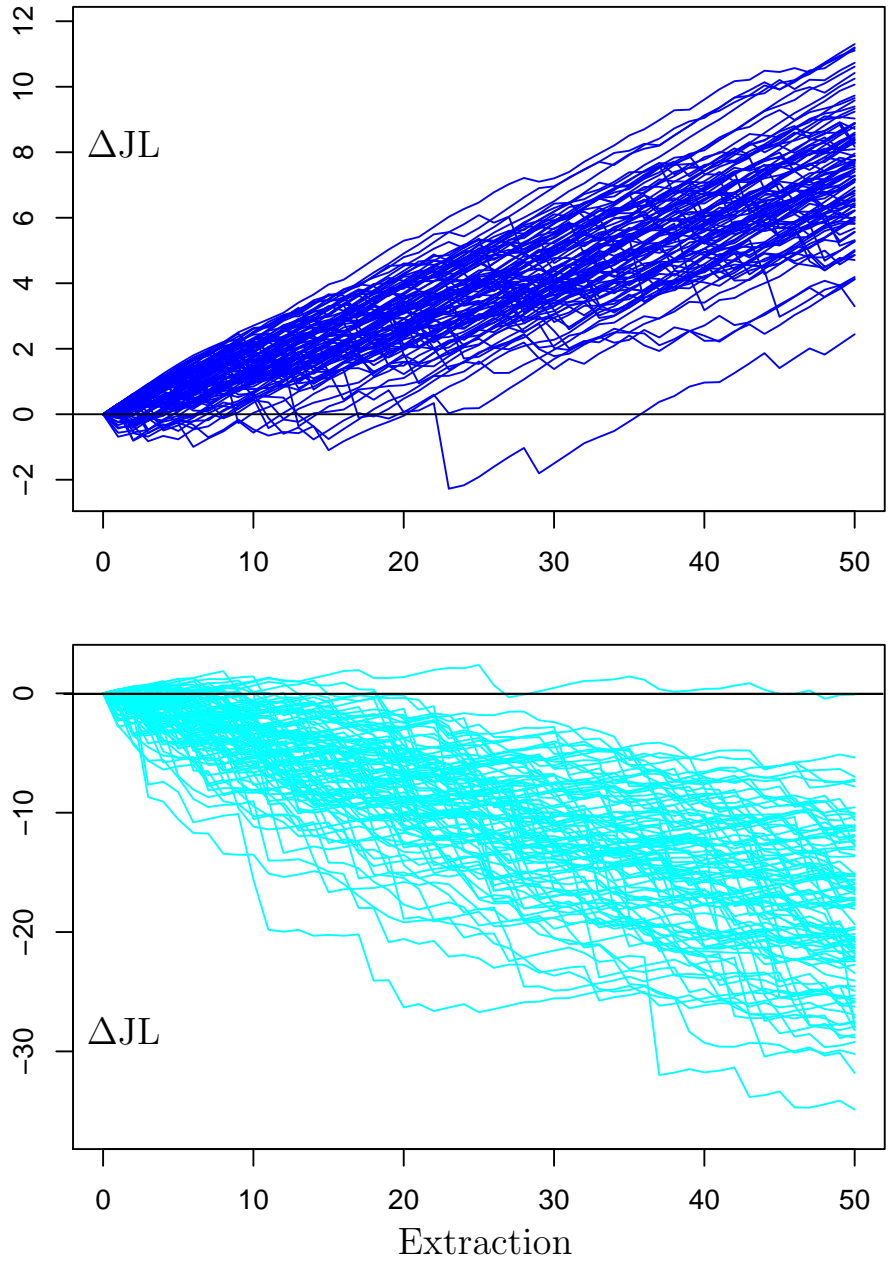,clip=,}
\caption{\small \sf Combined weights of evidence in simulated experiments.
The above (blue) combined JL sequences have been obtained by
the generator $H_1$, as it can be recognized because they
tend to large positive values as the number of extractions increases.
The below one are generated by $H_2$.}
\label{fig:simulazioni}
\end{figure}

Figure \ref{fig:simulazioni} shows 200 inferential stories, 
half per generator. We see that, in general, we get practically
sure of the model after a couple of dozens of extractions.
But there are also cases in which we need to wait longer
before we can feel enough sure on one hypothesis. 

It is interesting to remark that the leaning
in favor of each hypothesis grows, {\it in average}, linearly with
the number of extractions. That is, a little piece of evidence,
which is in average positive for $H_1$ and negative for $H_2$,
is added after each extraction. However, around the 
average trend, there is a large varieties of individual 
inferential histories. They all start at $\Delta\mbox{JL}=0$ for 
$n=0$, but in practice there are no two identical `trajectories'.
All together they form a kind of `fuzzy band',
whose `effective width' grows also with the number of extractions,
{\it but not linearly}. The widths grows as the square root 
of $n$.\footnote{We can evaluate the {\it prevision}
(`expected value')
of the variation of leaning at each random extraction for each hypotheses, 
calculated as the average value of $\Delta\mbox{JL}_{1,2}(H_i)$. 
We can also evaluate the {\it uncertainty of prevision},
quantified by the standard deviation. We get for the two hypotheses
\begin{eqnarray*}
\left\{\begin{array}{rcl} \mbox{E}[\Delta\mbox{JL}_{1,2}(H_1)] &=& 0.15 \\
                        \sigma[\Delta\mbox{JL}_{1,2}(H_1)] &= &0.24 
\\
u_R[\Delta\mbox{JL}_{1,2}(H_1)]  &=& 1.6 
%\frac{\sigma[\mbox{JL}_{1,2}(H_1)]}{|\mbox{E}[\mbox{JL}_{1,2}(H_1)]|} = 1.6
       \end{array}\right.
\hspace{0.5cm}&  &\hspace{0.5cm}
\left\{\begin{array}{rcl} \mbox{E}[\Delta\mbox{JL}_{1,2}(H_2)]&=&-0.38 \\
                        \sigma[\Delta\mbox{JL}_{1,2}(H_2)]&=&0.97 
\\
u_R[\Delta\mbox{JL}_{1,2}(H_2)] &=&2.6
%\frac{\sigma[\mbox{JL}_{1,2}(H_2)]}{|\mbox{E}[\mbox{JL}_{1,2}(H_2)]|} = 2.6
       \end{array}\right.
\end{eqnarray*} 
where also the {\it relative uncertainty} $u_R$ has been reported,
defined as the uncertainty divided by the absolute 
value of the prevision.
The fact that the uncertainties are relatively 
large tells clearly
that we {\it do not expect} that 
a single extraction will be sufficient to convince us of 
either  model.
But this does not mean we cannot take the decision
because the number of extraction {\it has been} too small.
If a very large fluctuation provides a $\Delta$JL of $-5$ 
(the table in this section shows that this is not very rare),
we have already got a very strong evidence in favor of $H_2$. 
Repeating what has been told several time, what matters
is the cumulated judgement leaning. It is irrelevant
if a  JL of $-5$ comes from ten individual pieces of evidence,
only from a single one, or partially from evidence
and partially from prior judgement.\\
When we plan to make $n$ extractions from a generator, 
probability theory allows
us to calculate expected value and uncertainty of $\mbox{JL}_{1,2}(n)$:
\begin{eqnarray*}
\mbox{E}[\Delta\mbox{JL}_{1,2}(n,H_i)] &=& n \times \mbox{E}[\Delta\mbox{JL}_{1,2}(H_i)] \\
\sigma[\Delta\mbox{JL}_{1,2}(n,H_i)] &=& \sqrt{n} \times \sigma[\Delta\mbox{JL}_{1,2}(H_i)]\\
u_R[\Delta\mbox{JL}_{1,2}(n,H_i)] &=& \frac{1}{\sqrt{n}}\times u_R[\Delta\mbox{JL}_{1,2}(H_i)]\,.
\end{eqnarray*}
In particular, for $n=50$ we get 
$\Delta\mbox{JL}_{1,2}(H_1) = 7.5\pm 1.7$ ($u_R=22\%$)
and  $\Delta\mbox{JL}_{1,2}(H_2) = -19\pm 7$ ($u_R=37\%$), that 
explain the gross feature of the bands in
 figure \ref{fig:simulazioni}.
}
This is the reason why, as $n$ increases, the bands tend 
to move away from the line $\mbox{JL}=0$. Nevertheless, individual
trajectories can exhibit very 
`irregular'\footnote{I find the issue of 
`statistical regularities' to be often misunderstood. 
For example, the trajectories in figure \ref{fig:simulazioni}
that do not follow the general trend {\it are not exceptions},
being generated by the same rules that produces all of them.} 
behaviors as we can  also see  in figure
\ref{fig:simulazioni}.
 
%\newpage
\section{Likelihood
%, likelihood principle
and maximum likelihood methods}
Some comments on {\it likelihood} are also in order,
because the reader might have heard this term and might
wonder if and how it fits in the scheme of reasoning expounded here.

One of the problems with this term is that it tends to have several
meanings, and then to create misunderstandings. 
In plane English `likelihood' is ``1. the condition of being likely or probable; probability'', or ``2. something that is 
probable''\footnote{See e.g. 
\url{http://www.thefreedictionary.com/likelihood}.};
but also ``3. (Mathematics \& Measurements / Statistics) the probability of a given sample being randomly drawn regarded as a function of the parameters of the population''.

Technically, with reference to the example of the previous appendix,
the likelihood is simply $P(x_E\,|\,H_i,I)$, where $x_E$ is fixed 
(the observation) and $H_i$ is the `parameter'. Then it can take two values, 
$P(x_E\,|\,H_1,I) = 3.68\times 10^{-8}$ and 
$P(x_E\,|\,H_2,I) = 1.99\times 10^{-8}$.

If, instead of 
only two models we had a continuity of models, 
for example
the family of
all Gaussian distributions characterized by central value $\mu$
and `effective width' (standard deviation) $\sigma$, our likelihood
would be $P(x_E\,|\,\mu,\sigma,I)$, i.e. 
\begin{eqnarray}
{\cal L}(\mu,\sigma\,;\,x_E) &=& P(x_E\,|\,\mu,\sigma,I)\,,
\end{eqnarray}
written in this way to remember that: 1) a likelihood is a function
of the model parameters and not of the data; 
2) ${\cal L}(\mu,\sigma\,;\,x_E)$ {\it is not} a probability
(or a probability density function) of $\mu$ and $\sigma$.
Anyway, for the rest of the discussion we stick to 
the very simple likelihood based on the two Gaussians.
That is, instead of a double infinity of possibilities, 
our space of parameters is made only of two points, $\{\mu_1=0,\sigma_1=1\}$
and  $\{\mu_1=0.4,\sigma_2=2\}$. Thus the situation gets simpler,
although the main conceptual issues remain substantially the same.

In principle there is nothing bad to give a special name 
to this function of the parameters. But, frankly, I had preferred
statistics gurus named it after their dog
or their lover, rather than call it
`likelihood'.\footnote{Note added: I have just learned, 
while making the short 
research on the use of the logarithmic updating of the odds
presented in Appendix E, that ``the term [likelihood] 
was introduced by R. A. Fisher with the object of 
{\it avoiding} the use of Bayes' theorem''~\cite{Good}.
}
The problem is that it is very frequent to
hear students, teachers and researcher
 explaining that the `likelihood' tells 
``how likely the parameters are'' (this is {\it the probability
of the parameters! not the `likelihood'}). Or they would say,
with reference to our example, 
``it is the probability that $x_E$ comes from $H_i$'' 
(again, this expression  would be the probability of $H_i$ 
given $x_E$, and not the probability of $x_E$ given the 
models!) Imagine if we have only $H_1$ in the game:
$x_E$ comes with certainty from $H_1$,
although $H_1$ does not yield with certainty $x_E$.\footnote{
As further example, you might look at
\url{http://en.wikipedia.org/wiki/Likelihood_principle}, where
it is stated (January 28, 2010, 15:40)
that a likelihood
``gives a measure of how `likely' any particular value of $\theta$ is''
(note the quote mark of `likely', as in the example of footnote
\ref{fn:james}).
But, fortunately we
 find in \url{http://en.wikipedia.org/wiki/Likelihood_function} that
``This is not the same as the probability that those parameters
 are the right ones, given the observed sample. 
Attempting to interpret the likelihood of a hypothesis 
given observed evidence as the probability of the hypothesis 
is a common error, with potentially disastrous real-world 
consequences in medicine, engineering or jurisprudence. 
See prosecutor's fallacy[*] for an example of this.''
([*] see 
\url{http://en.wikipedia.org/wiki/Prosecutor\%27s_fallacy}.)
\\
Now you might understand why I am particular upset 
with the name likelihood. 
\label{fn_wiki_likelihood}
} 

Several methods in `conventional statistics'
use somehow the likelihood to decide which model
or which set of parameters describes at best the data.
Some even use the likelihood ratio (our Bayes factor),
or even the logarithm of it (something equal or proportional,
depending on the base, to the weight of evidence we have indicated
here by JL). 
The most famous method of the series
is the {\it maximum likelihood principle}.
As it is easy to guess from its name, it states
that the {\it best estimates} of the parameters
are those which maximize the likelihood.

All that {\it seems} reasonable and in agreement 
with what it has been expounded here, but it is not
quite so. First, for those who support this approach,
likelihoods are not just a part of the inferential
tool, they are everything. Priors are completely
neglected, more or less because of the objections
in footnote \ref{fn:nopriors}. This can be acceptable, 
if the evidence is overwhelming, but this is not always the case.
Unfortunately, as it is now easy to understand, neglecting
priors  is mathematically
equivalent to consider the alternative hypotheses equally likely!
As a consequence of this statistics miseducation 
(most statistics courses in the universities all around the world
only teach `conventional statistics' and never, little, or
badly probabilistic inference)
is that too many unsuspectable people 
fail in solving the AIDS problem of appendix B, 
or confuse the likelihood with the probability of the 
hypothesis, resulting in misleading scientific claims
 (see also footnote \ref{fn_wiki_likelihood} and Ref.~\cite{BR}). 

The second difference is that, since ``there are no priors'',
the result cannot have a probabilistic meaning, as
it is
openly recognized by the promoters of this method,
who, in fact, do not admit we can talk about probabilities of causes
(but most practitioners seem not to be aware of this 
`little philosophical detail', also because frequentistic 
gurus, having difficulties to explain what is the meaning
of their methods, they say they are `probabilities', 
but in quote marks!\footnote{
For example,
we read in Ref.~\cite{JamesRoos} 
(the authors are influential supporters
of the use frequentistic methods in the particle physics community):
\begin{quote}
{\sl 
When the result of a measurement of a physics quantity is published as 
$R=R_0\pm\sigma_0$ without further explanation, it simply implied
that R is a Gaussian-distributed measurement with mean $R_0$
and variance $\sigma_0^2$. 
This allows to calculate various confidence
intervals of given ``probability'', i.e. the ``probability''
P that the true value of $R$ is within a given interval.
}
\end{quote} 
(Quote marks are original and nowhere in the paper is
explained why probability is in quote marks!)\\
The following Good's words about frequentistic
{\it confidence intervals} (e.g. `$R=R_0\pm\sigma_0$' of the 
previous citation) and ``probability''
might be very enlighting (and perhaps shocking, 
if you always thought they meant something like 
`how much one is confident in something'):
\begin{quote}
{\sl
Now suppose that the functions $\underline{c}(E)$ 
and  $\overline{c}(E)$ are selected so that 
$[\overline{c}(E),\overline{c}(E)]$ is a confidence interval with 
coefficient $\alpha$, where $\alpha$ is near to 1. Let us assume 
that the following instructions are issued to all statisticians.

``Carry out your experiment, calculate the confidence interval,
and {\it state} that $c$ belong to this interval.
If you are asked whether you `believe' that
$c$ belongs to the confidence interval you must refuse to answer.
In the long run your assertions, if independent of each other, will
be right in approximately a proportion $\alpha$ of cases.''
(Cf. Neyman (1941), 132-3)
}\cite{Good}
\end{quote}
[Neyman (1941) stands for J. Neyman's ``Fiducial argument
and the theory of confidence intervals'', {\it Biometrica}, 
{\bf 32}, 128-150.]\\
(For comments about what is in my opinion
a ``kind of condensate of frequentistic nonsense'',
see Ref.~\cite{BR}, in particular section 10.7
on {\it frequentistic coverage}. You might
get a feeling of what happens
taking Neyman's prescriptions literally
playing with
the `the ultimate confidence intervals calculator' 
available in
\url{http://www.roma1.infn.it/~dagos/ci_calc.html}.)  
\label{fn:james}}).
As a consequence, the resulting `error analysis',
that in human terms means to assign different 
beliefs to different values of the parameters,
is cumbersome. In practice the results are reasonable only
if the possible values of the parameters are
initially equally likely and the `likelihood function' has
a `kind shape' (for more details see chapters 1 and 12
of Ref.~\cite{BR}).
\\ \break \mbox{}
%\newpage
\section{Evidences mediated by a testimony}
In most cases (and practically always in courts)
pieces of evidence are not acquired directly by the person
who has to form his mind about the plausibility of a hypothesis.
They are usually accounted by an intermediate person, or
by a chain of individuals. 
Let us call $E_T$ the report of the 
evidence $E$ provided in a {\it testimony}.
The inference becomes now $P(H_i\,|\,E_T,I)$, generally
different from $P(H_i\,|\,E,I)$.

In order to apply Bayes' theorem in one of its form
we need first to evaluate $P(E_T\,|\,H_i,I)$. 
Probability theory teaches us how to get it 
[see Eq.~(\ref{eq:rul4}) in Appendix A]:
\begin{eqnarray}
P(E_T\,|\,H_i,I) &=&    P(E_T\,|\,E,I)\cdot  P(E\,|\,H_i,I)
                     +  P(E_T\,|\,\overline E,I)\cdot  P(\overline E\,|\,H_i,I)
\end{eqnarray}
($E_T$ could be due to a true evidence or to a fake one).
Three new ingredients enter the game:
\begin{itemize}
\item  $P(E_T\,|\,E,I)$, that is the probability of
       the evidence to be correctly reported as such.
\item  But the testimony could also be incorrect the other way around
       (it could be incorrectly reported, simply by mistake, 
       but also it could be a `fabricated evidence'),
       and therefore also $P(E_T\,|\,\overline E,I)$ is needed. 
       Note that the probabilities to lie could be in general
       asymmetric, i.e. $P(\overline E_T\,|\,E,I)\ne P(E_T\,|\,\overline E,I)$,
       as we have seen in the AIDS problem of Appendix F, in which 
       the response of the analysis models false witness well. 
\item  Finally, since  $P(\overline E\,|\,H_i,I)$ enters now directly,
       the `na\"\i ve' Bayes factor, only depending 
       on $P(E\,|\,H_i,I)$, is not longer enough. 
\end{itemize}
Taking our usual two hypotheses, $H_1=H=\mbox{`guilty'}$
and  $H_2=\overline{H}=\mbox{`innocent'}$, we get the following
Bayes factor based on the {\it testified evidence} $E_T$
(hereafter, in order to simplify the notation,
 we use the subscript `$H$' in odds and Bayes factors, 
instead of `$i,j$', 
to indicate that they are in favor of $H$ and against $\overline H$,
as we already did in the 
AIDS example of Appendix F):
\begin{eqnarray}
\tilde O_H(E_T,I) &=& \frac{ P(E_T\,|\,E,I)\cdot  P(E\,|\,H,I)
                     +  P(E_T\,|\,\overline E,I)\cdot  P(\overline E\,|\,H,I)
                        }
                        { P(E_T\,|\,E,I)\cdot  P(E\,|\,\overline H,I)
                     +  P(E_T\,|\,\overline E,I)\cdot  P(\overline E\,|\,\overline H,I)
                        }\,.\label{eq:BF_lie}
\end{eqnarray}
As expected, this formula is a bit more complicate that the 
Bayes factor calculated taking $E$ for granted, which is recovered
if the lie probabilities vanish
\begin{eqnarray} 
\tilde O_{H}(E_T,I) && \xrightarrow [\mbox{{\footnotesize 
             $ P(E_T\,|\,\overline E,I)\rightarrow 0$}}]{}
                                    {\ \ \ \tilde O_{H}(E,I)}\,,
\end{eqnarray}
i.e. only  when 
we are absolutely sure the witness does not 
err or lie reporting $E$
(but Peirce reminds us that ``absolute certainty, 
or an infinite chance, can never be attained by 
mortals''~\cite{Peirce}).

In order to single out the effects of the new ingredients, Eq.~(\ref{eq:BF_lie})
can be rewritten as\footnote{
Factorizing $P(E\,|\,H,I)$ and $P(E\,|\,\overline H,I)$ respectively
in the numerator and in the denominator, Eq.~(\ref{eq:BF_lie}) 
becomes
\begin{eqnarray*}
\tilde O_{H}(E_T,I)
                    &=& \tilde O_{H}(E,I)\times
                        \frac{1 +  \frac{P(E_T\,|\,\overline E,I)}{P(E_T\,|\,E,I)} \cdot 
                              \frac{P(\overline E\,|\,H,I)}{P(E\,|\,H,I)}
                             }
                             {1 +  \frac{P(E_T\,|\,\overline E,I)}{P(E_T\,|\,E,I)} \cdot 
                               \frac{P(\overline E\,|\,\overline H,I)}{P(E\,|\,\overline H,I)}
                             }\,.
\end{eqnarray*}
Then $P(E_T\,|\,\overline E,I)/P(E_T\,|\,E,I)$ can be indicated as $\lambda(I)$, 
$P(\overline E\,|\,H_i,I)$ is equal to $1-P(E\,|\,H_i,I)$ and, finally,
$P(E\,|\,\overline H,I)$ can be written as  $P(E\,|\,H,I)/\tilde O_{H}(E,I)$.
}
\begin{eqnarray}
\tilde O_{H}(E_T,I)
                    &=& \tilde O_{H}(E,I)\times
                        \frac{1 +  \lambda(I) \cdot \left[
                               \frac{1}{P(E\,|\,H,I)} -1
                               \right]}
                             {1 +  \lambda(I) \cdot \left[
                                \frac{ \tilde O_{H}(E,I)}{P(E\,|\,H,I)} -1
                               \right]}
\,, \label{eq:weight_ET}
\end{eqnarray}
where 
\begin{eqnarray}
\lambda(I) &=& \frac{P(E_T\,|\,\overline E,I)}{P(E_T\,|\,E,I)}\,,
\end{eqnarray}
under the {\it condition}\footnote{Otherwise, obviously  
$\tilde O_{H}(E,I)$ cannot be factorized. 
The effective odds $\tilde O_{H}(E_T,I)$ can 
however be written in the following convenient forms
\begin{eqnarray*}
\left.\tilde O_{H}(E_T,I)\right|_{P(E\,|\,H,I) = 0}
&=& \frac{1}{P(\overline E\,|\,\overline H)+P(E\,|\,\overline H)/\lambda}\\
\left.\tilde O_{H}(E_T,I)\right|_{P(E\,|\,\overline H,I) = 0}
&=& \lambda\,P(E\,|\,H)+P(\overline E\,|\,H)\,,
\end{eqnarray*}
although less interesting than Eq.~(\ref{eq:weight_ET}).
\label{fb:PEH=0}
} $P(E\,|\,H,I)> 0$ and
$P(E\,|\,\overline H,I)> 0$, i.e. 
{\it $\tilde O_{H}(E,I)$ positive and finite}.
The parameter
$\lambda(I)$, ratio of the {\it probability of fake evidence} and the 
{\it probability that the evidence is correctly accounted}, 
can be interpreted as a kind of {\it lie factor}.
Given the human roughly logarithmic sensibility to probability
ratios, it might be useful to define, in analogy to the JL,
\begin{eqnarray}
\mbox{J}\lambda(I) &=& \log_{10}[\lambda(I)]\,.
\end{eqnarray}
Let us make some instructive limits of Eq.~(\ref{eq:weight_ET}).
\begin{eqnarray}
\tilde O_{H}(E_T,I) &\xrightarrow [\mbox{{\footnotesize 
                        $\lambda(I) \rightarrow 0$}}]{}{}&\tilde O_{H}(E,I)\\
\tilde O_{H}(E_T,I) &\xrightarrow [\mbox{{\footnotesize 
                        $\lambda(I) \rightarrow 1$}}]{}{}& 1 \\
\tilde O_{H}(E_T,I) & \xrightarrow [\mbox{{\footnotesize 
             $P(E\,|\,H,I)\rightarrow 0$}}]{}{}& 1\\
\tilde O_{H}(E_T,I) & \xrightarrow [\mbox{{\footnotesize 
             $\tilde O_{H}(E,I)\rightarrow \infty $}}]{}{}& 
           \frac{P(E\,|\,H,I)}{\lambda(I)} + 1 - P(E\,|\,H,I)
\end{eqnarray}
As we have seen, the ideal case is recovered if the lie factor vanishes.
Instead, if it is equal to 1, i.e. $\mbox{J}\lambda(I)=0$, the reported
evidence becomes useless. The same happens if $P(E\,|\,H,I)$ vanishes 
[this implies that $P(E\,|\,\overline H,I)$ vanishes too, being
$P(\overline H,I)=P(E\,|\,\overline H,I)/\tilde O_{H}(E,I)$].

However, the most remarkable limit
is the last one.  It states 
that, even if $\tilde O_{H}(E,I)$ is very high, 
the effective Bayes factor cannot exceed the inverse of the lie factor:
\begin{eqnarray}
\tilde O_{H}(E_T,I) & \le & \frac{P(E\,|\,H,I)}{\lambda(I)} \le \frac{1}{\lambda}\hspace{0.5cm}\mbox{[if $\tilde O_{H}(E,I)\rightarrow\infty$]}\,,
\end{eqnarray}
or, using logarithmic quantities
 \begin{eqnarray}
\Delta\mbox{JL}(E_T,I) & \le & - \mbox{J}\lambda + \log_{10}{P(E\,|\,H,I)} \le
  - \mbox{J}\lambda
\hspace{0.5cm}\mbox{[if $\Delta$JL$(E,I)\rightarrow\infty$]} \,.
\end{eqnarray}
At this point some numerical examples are in order (and those
who claim they can form their mind on pure intuition
get all my admiration\ldots {\it if} they really can).
Let us imagine that $E$ would ideally provide a weight of evidence 
of 6 [i.e. $\Delta\mbox{JL}_H(E,I)=6$]. We can study, with the help 
of table \ref{tab:JL_ET},
\begin{table}[!t]
\begin{center}
\begin{tabular}{|c|c|c|c|c|c|c|c|c|c|}
\multicolumn{10}{c}{}\\ 
\hline
\multicolumn{10}{|c|}{$\Delta\mbox{JL}_H(E,I)=6$} \\
\hline
$\mbox{J}\lambda(I)$
                  && \multicolumn{8}{|c|}{$\mbox{JL}_H(E_T,I)$} \\
\hline
                 &$\mbox{JL}_{E}(H,I)$:& 10 & $3$ & $2$ & $1$ & 
                                        $0$ & $-1$ & $-3$ & $-10$
\\
\hline 
$\rightarrow -\infty$ && 6.00 & 6.00 & 6.00 & 6.00 & 6.00 & 6.00 & 6.00 & 6.00 \\
$-8$                    && 6.00 & 6.00 & 6.00 & 6.00 & 5.99 & 5.95 & 4.96 & $4\times 10^{-3}$ \\
$-7$                    && 5.96 & 5.96 & 5.96 & 5.95 & 5.92 & 5.68 & 4.00 & $4\times 10^{-4}$\\
$-6$                    && 5.70 & 5.70 & 5.70 & 5.68 & 5.52 & 4.92 & 3.00 & $4\times 10^{-5}$\\
$-5$                    && 4.96 & 4.96 & 4.95 & 4.92 & 4.68 & 3.95 & 2.00 & $4\times 10^{-6}$\\
$-4$                    && 4.00 & 4.00 & 3.99 & 3.95 & 3.70 & 2.96 & 1.04 & $4\times 10^{-7}$\\
$-3$                    && 3.00 & 3.00 & 3.00 & 2.96 & 2.70 & 1.96 & 0.30 & $4\times 10^{-8}$\\
$-2$                    && 2.00 & 2.00 & 2.00 & 1.95 & 1.70 & 1.00 & 0.04 & $4\times 10^{-9}$\\
$-1$                    && 1.00 & 1.00 & 1.00 & 0.96 & 0.74 & 0.26 & 0.004& $4\times 10^{-10}$\\
%$-0.8$                  && 0.80 & 0.80 & 0.80 & 0.77 & 0.56 & 0.17 & 0.002& $2\times 10^{-10}$\\  
%$-0.6$                  && 0.60 & 0.60 & 0.60 & 0.57 & 0.40 & 0.10 & 0.001& $1\times 10^{-10}$\\  
%$-0.4$                  && 0.40 & 0.40 & 0.40 & 0.38 & 0.24 & 0.06 & 0.0007& $7\times 10^{-11}$\\
%$-0.2$                  && 0.20 & 0.20 & 0.20 & 0.19 & 0.11 & 0.02 & 0.0003& $3\times 10^{-11}$\\
\ 0                   && \ 0  & \ 0  & \ 0  & \ 0    & \ 0    & \ 0   &\  0 &\ 0    \\
\hline
 \multicolumn{10}{c}{}\\ 
\hline
 \multicolumn{10}{|c|}{$\Delta\mbox{JL}_H(E,I)=3$} \\
\hline
$\mbox{J}\lambda(I)$
                  && \multicolumn{8}{|c|}{$\mbox{JL}_H(E_T,I)$} \\
 
\hline
                 &$\mbox{JL}_{E}(H,I)$:& 10 & $3$ & $2$ & $1$ & 
                                        $0$ & $-1$ & $-3$ & $-10$
\\
\hline 
$\rightarrow -\infty$ && 3.00 & 3.00 & 3.00 & 3.00 & 3.00 & 3.00 & 3.00 & 3.00 \\
%$-6$                    && 3.00 & 3.00 & 3.00 & 3.00 & 3.00 & 3.00 & 2.70 & $4\times 10^{-5}$ \\
$-5$                    && 3.00 & 3.00 & 3.00 & 3.00 & 2.99 & 2.95 & 1.96 & $4\times 10^{-6}$ \\
$-4$                    && 2.96 & 2.96 & 2.96 & 2.95 & 2.92 & 2.68 & 1.04 & $4\times 10^{-7}$ \\
$-3$                    && 2.70 & 2.70 & 2.70 & 2.68 & 2.52 & 1.93 & 0.30 & $4\times 10^{-8}$ \\
$-2$                    && 1.96 & 1.96 & 1.96 & 1.92 & 1.68 & 1.00 & 0.04 & $4\times 10^{-9}$\\
$-1$                    && 1.00 & 1.00 & 0.99 & 0.96 & 0.74 & 0.26 & 0.004 & $4\times 10^{-10}$\\
\ 0                   && \ 0  & \ 0  & \ 0  & \ 0    & \ 0    & \ 0   &\  0 &\ 0    \\
\hline\multicolumn{10}{c}{}\\ 
\hline
 \multicolumn{10}{|c|}{$\Delta\mbox{JL}_H(E,I)=1$} \\
\hline
$\mbox{J}\lambda(I)$
                  && \multicolumn{8}{|c|}{$\mbox{JL}_H(E_T,I)$} \\
 
\hline
                 &$\mbox{JL}_{E}(H,I)$:& 10 & $3$ & $2$ & $1$ & 
                                        $0$ & $-1$ & $-3$ & $-10$
\\
\hline 
$\rightarrow -\infty$ && 1.00 & 1.00 & 1.00 & 1.00 & 1.00 & 1.00 & 1.00 & 1.00 \\
$-3$                    && 1.00 & 1.00 & 1.00 & 1.00 & 0.99 & 0.96 & 0.26 & $4\times 10^{-8}$ \\
$-2$                    && 0.96 & 0.96 & 0.96 & 0.96 & 0.93 & 0.72 & 0.04 & $4\times 10^{-9}$\\
$-1$                    && 0.72 & 0.72 & 0.72 & 0.70 & 0.58 & 0.23 & 0.003 & $4\times 10^{-10}$\\
$-0.5$                  && 0.41 & 0.41 & 0.41 & 0.39 & 0.27 & 0.07 &  $8\times 10^{-4}$  & $4\times 10^{-10}$\\
\ 0                   && \ 0  & \ 0  & \ 0  & \ 0    & \ 0    & \ 0   &\  0 &\ 0    \\
\hline
\end{tabular}
\caption{\small \sf Dependence 
of the judgement leaning due to a reported evidence 
[$\Delta\mbox{JL}_H(E_T,I)$] for $\Delta\mbox{JL}_H(E,I)=6$, 3 and 1
as a function the other ingredients of the inference
(see text). Note the upper 
limit of $\Delta\mbox{JL}_H(E_T,I)$ to 
 $-\mbox{J}\lambda$, if this value is $\le \Delta\mbox{JL}_H(E,I)$. 
}
\label{tab:JL_ET}
\end{center}
\end{table}
 how the {\it weight of the reported evidence} $\Delta\mbox{JL}_H(E_T,I)$
depends on the other beliefs [in this table logarithmic
quantities have been used throughout, therefore 
$\mbox{JL}_E(H,I)$ is the base ten logarithm of the odds in favor
of $E$ given the hypothesis $H$; the table provides, for comparisons,
also  $\Delta\mbox{JL}_H(E_T,I)$ from  
$\Delta\mbox{JL}_H(E,I)$ equal to 3 and 1]. 

The table exhibits the limit behaviors we have seen analytically. 
In particular,
if we fully trust the report, i.e. $\mbox{J}\lambda(I)=-\infty$,
then $\Delta\mbox{JL}_H(E_T,I)$ is exactly equal to 
$\Delta\mbox{JL}_{H}(E,I)$,
as we already know. But as soon as the absolute value of the lie factor
is close to $\mbox{JL}_H(E,I)$, there is a sizeable effect.
The upper bound can be the be rewritten as
\begin{eqnarray}
\tilde O_{H}(E_T,I) & \le & \mbox{min}\,
[\tilde O_{H}(E,I),\, \frac{1}{\lambda}]\,,\\
\mbox{or}\hspace{3cm}\mbox{} && \mbox{} \nonumber\\
   \Delta\mbox{JL}_H(E_T,I) & \le & 
       \mbox{min}\,[\Delta\mbox{JL}_H(E,I),\,-\mbox{J}\lambda(I)] \,,
\end{eqnarray}
a relation valid in the region of interest when thinking 
about an evidence in favor of $H$, i.e. $\Delta\mbox{JL}_H(E,I) > 0$
and $\mbox{J}\lambda(I) < 0$. 

This upper bound is very interesting. Since minimum conceivable values 
of $\mbox{J}\lambda(I)$ for human beings can be of the order of
$-6$ (to perhaps  $\approx -8$ or  $\approx -9$,
but in many practical applications $-2$ or  $-3$ can already be very generous!), 
in practice the effective weights of evidence cannot exceed
values of about $+6$ (I have no strong opinion on the exact
value of this limit, my main point is that {\it you consider 
there might be
such a practical human limit}.)

This observation has an important consequence in the combination of evidences, 
as anticipated at the end of section \ref{ss:Agatha}. 
Should we give more consideration to
a single strong piece of evidence, virtually weighing 
$\Delta\mbox{JL}(E)=10$, or 10 independent weaker evidences, 
each having a $\Delta$JL of 1? As it was said, in the ideal case
they yield the same global leaning factor. 
But as soon as human fallacy (or conspiracy) is taken into account,
and we remember that our belief is based on $E_T$ and not on $E$,
then we realize that 
$\Delta\mbox{JL}(E_T)=10$ is well above the range 
of JL that we can reasonably conceive. 
Instead the weaker pieces of evidence are little affected
by this doubt and when they sum up together, they really 
can provide a $\Delta$JL of about 10. 
%\\ \break \mbox{}

%\newpage \mbox{}
%\newpage 
\section{A simple Bayesian network}
\begin{figure}[!t]
\centering\epsfig{file=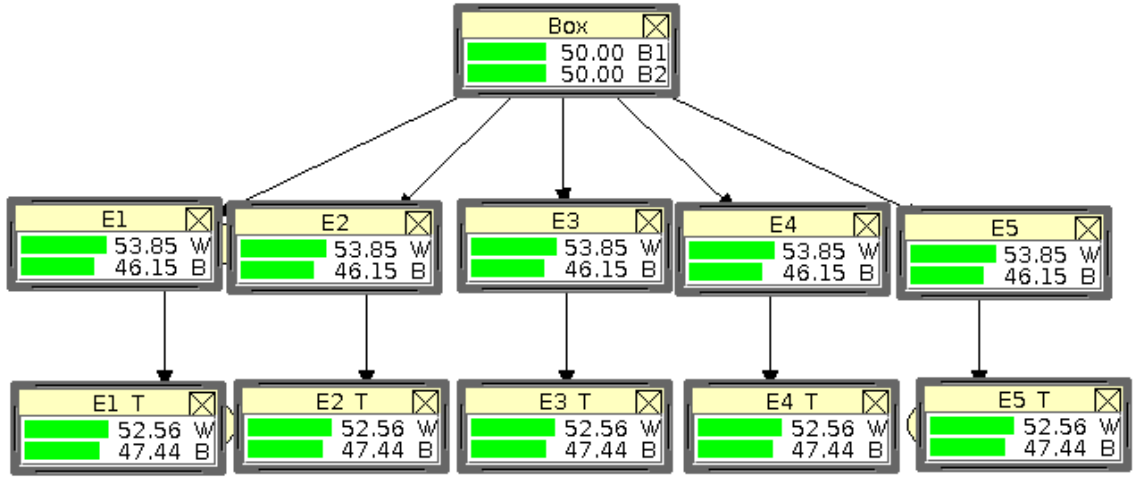,clip=,}
\caption{\small \sf Same belief network of figure \ref{fig:bn}. 
This representation shows the `monitors' giving the initial probabilities
of all states of the variables. {\bf If you like to test your intuition},
try the guess how {\it all} probabilities change when the information 
changes in the following order: a) witness 1 says white ($E_1T=W$); 
b) witness 2 also reports white ($E_2T=W$); c) then witness 3
claims, contrary to the previous two, to have observed
black ($E_3T=B$);
 c) finally we directly observe the result of the fourth extraction, resulting black 
($E_4=B$). The {\bf solutions} are in figures \ref{fig:bn_mon_ET1_w} to 
\ref{fig:bn_mon_ET_ww_ET_b_Eb}.
}
\label{fig:bn_mon_0}
\end{figure}

\begin{figure}
\centering\epsfig{file=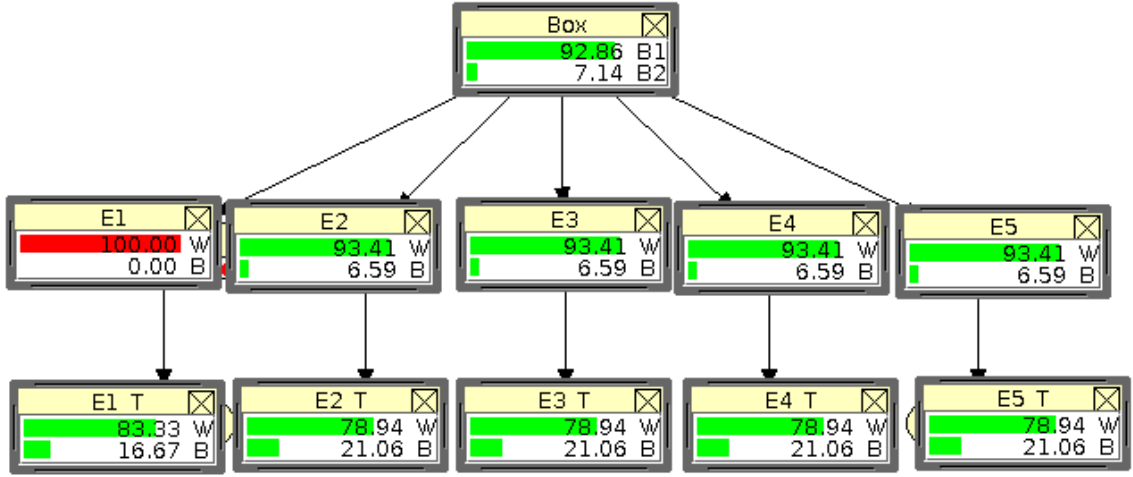,clip=,}
\caption{\small \sf  Status of the network after the observation 
of a white ball.}
\label{fig:bn_mon_E1_w}
\end{figure}

Let us go back to our toy model of section \ref{sec:1in13}
and let us complicate it just a little bit, adding the possibility 
of incorrect testimony (but we also simplify it
using uniform priors, so that we can focus
on the effect of the {\it uncertain evidence}). 
For example, imagine you do not see
directly the color of the ball, but this is reported to you by
a collaborator, who, however, might not tell you always the truth.
We can model the possibility of a lie in following way:
after each extraction he tosses a die and reports the true 
color only if the die gives a number smaller than 6. 
Using the formalism 
of Appendix I, we have
\begin{eqnarray}
P(E_T\,|\,E,I) &=& 5/6 \\
P(E_T\,|\,\overline E,I) &=& 1/6\,, \\
\mbox{i.e.}\hspace{3cm} && \hspace{3cm}\nonumber \\
\lambda(I) &=& 1/5.
\end{eqnarray}
The resulting {\it belief network},\footnote{In complex situations
an effects might have several (con-)causes; or
an effect can be itself a cause of other effects;
and so on. As it can be easily imagined, 
causes and effects can be represented by a graph, 
as that of figure   \ref{fig:bn}.
Since the connections between the {\it nodes}
of the resulting {\it network} have usually the meaning
of probabilistic links (but also deterministic relations
can be included), this graph is called a {\it belief network}.
Moreover, since Bayes' theorem is used to update the probabilities 
of the possible {\it states} of the nodes (the node `Box', 
with reference to our toy model, has states $B_1$ and $B_2$;
the node `Ball' has states $W$ and $B$), they are also
called {\it Bayesian networks}. 
For more info, as well as tutorials and demos 
of  powerful packages  having 
also a friendly graphical user interface, 
I recommend visiting Hugin~\cite{Hugin} and
Netica~\cite{Netica} web sites.
(My preference for Hugin is mainly due to the fact that 
it is multi-platform and runs nicely under Linux.)
For a book 
introducing Bayesian networks in forensics, Ref.~\cite{Taroni}
is recommended. For a monumental 
probabilistic network on the `{\bf case that will never end}', 
see Ref.~\cite{SaccoVanzetti} 
(if you like classic thrillers, 
the recent paper of the same author might be of your interest
\cite{JBK}).}
relative to five extractions and
to the corresponding five reports is shown in 
figure \ref{fig:bn},
redrawn in a different way  
in figure \ref{fig:bn_mon_0}. 
In this diagram 
the {\it nodes} are represented by `monitors' that
provide the probability of each {\it state} of the {\it variable}.
The green bars mean that we are in condition of uncertainty
with respect to all states of all variable. 
Let us describe the several nodes:
\begin{itemize}
\item
Initial box compositions have probability 50\% each,
that was our assumption. 
\item
The probability of white and black are the same for all
extractions, with white a bit more probable than black
(14/26 versus 12/26, that is 53.85\% versus 46.15\%).
\item
There is also higher probability that the `witness' reports
white, rather than black, but the difference is attenuated by
the `lie factors'.\footnote{Note that there are in general two lie factors,
one for $E$ and one for $\overline E$. For simplicity we assume here
they have the same value.} 
In fact, calling $W_T$ and $B_T$ the reported
colors we have
\begin{eqnarray}
P(W_T\,|\,I) &=& P(W_T\,|\,W,I)\cdot P(W\,|\,I) + 
                 P(W_T\,|\,B,I)\cdot P(B\,|\,I) \\
P(B_T\,|\,I) &=& P(B_T\,|\,W,I)\cdot P(W\,|\,I) + 
                 P(B_T\,|\,B,I)\cdot P(B\,|\,I)\,.
\end{eqnarray}
\end{itemize}
Let us now see what happens if we {\it observe white}
(red bar in figure \ref{fig:bn_mon_E1_w}).
All probabilities of the network have been updated
(Hugin~\cite{Hugin} has  nicely done the job 
for us\footnote{The Hugin file
can be found in 
\url{http://www.roma1.infn.it/~dagos/prob+stat.html\#Columbo}.}). 
We recognize the 93\% of box $B_1$, that we already know. 
We also see that the increased belief on this box makes
us more confident to observe white balls in the following
extractions (after re-introduction). 

\begin{figure}
\centering\epsfig{file=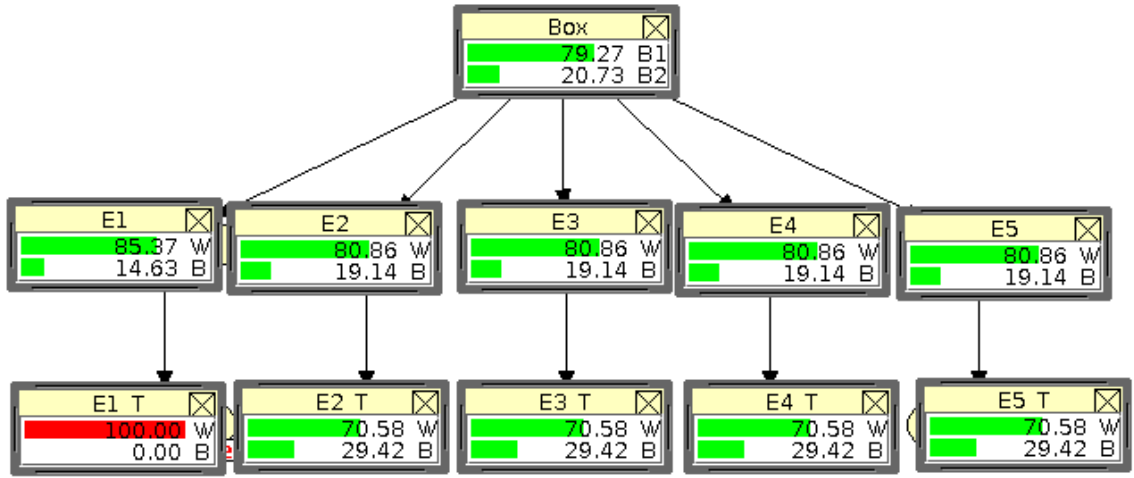,clip=,}
\caption{\small \sf Status of the network after the report
of a white ball (compare with figure \ref{fig:bn_mon_E1_w}).}
\label{fig:bn_mon_ET1_w}
\end{figure}

More interesting is the case in which our inference
is based on the reported color (figure \ref{fig:bn_mon_ET1_w}).
The fact that the witness could lie reduces,
with respect to the previous case, our confidence
on $B_1$  and on white balls
in future extractions. As an exercise on what we have learned 
in appendix H, we can evaluate the `effective' Bayes factor
$\tilde O_{B_1}(W_T,I)$ that takes into account the testimony.
Applying Eq.~(\ref{eq:weight_ET}) we get
\begin{eqnarray}
\tilde O_{B_1}(W_T,I) &=& \tilde O_{B_1}(W,I) \times 
               \frac{1 +  \lambda(I) \cdot \left[
                               \frac{1}{P(W\,|\,B_1,I)} -1
                               \right]}
                             {1 +  \lambda(I) \cdot \left[
                                \frac{ \tilde O_{H}(W,I)}{P(W\,|\,H,I)} -1
                               \right]} \\
&=& 13\times \frac{5}{17} = 3.82\,,
\end{eqnarray}
or $\Delta\mbox{JL}_{B_1}(W_T,I)=0.58$, about a factor of two
smaller than $\Delta\mbox{JL}_{B_1}(W,I)$, that was 1.1
(this mean we need two pieces of evidence of this kind
to recover the loss of information due to the testimony).  

The network gives us also the probability that the witness
has really told us the truth, i.e. $P(W\,|\,W_T,I)$, that is
{\it different} from $P(W_T\,|\,W,I)$, the reason being that
white was initially a bit more probable than black. 

\begin{figure}
\centering\epsfig{file=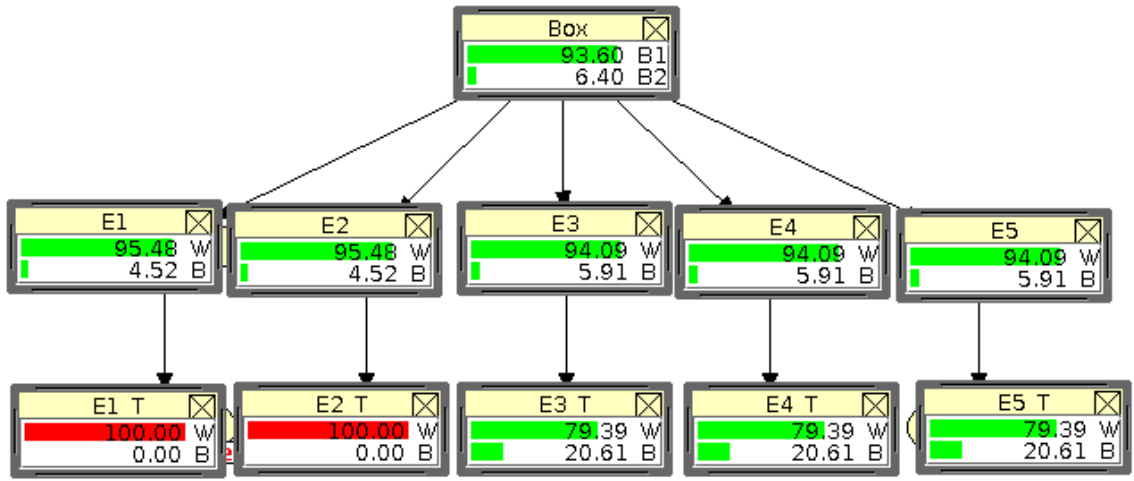,clip=,}
\caption{\small \sf Network of figure \ref{fig:bn_mon_ET1_w} updated
by a second testimony in favor of white.}
\label{fig:bn_mon_ET_ww}
\end{figure}

Let us see now what happens if we get two concording testimonies 
(figure \ref{fig:bn_mon_ET_ww}).
As expected, the probability of $B_1$ increases and becomes
closer to the case of a direct observation of white. 
As usual, also the probabilities of future white balls increase.

The most interesting thing that comes from the
result of the network is how the probabilities 
that the two witness lie change. First we see that they are the same, 
about 95\%, as expected for symmetry. But the surprise is that
the probability the the first witness said the truth has increased,
passing from 85\% to 95\%. We can justify the variation
because, in qualitative agreement with intuition, if we have 
concordant witnesses, {\it we tend to believe to each of them more
than what we believed individually}. Once again, the result is, 
perhaps after an initial surprise, in qualitative agreement with
intuition. The important point is that intuition is unable to
get quantitative estimates. Again, the message is that,
once we agree on the basic assumption and we check, whenever it
is possible, that the results are reasonable, it is better to rely on 
automatic computation of beliefs.
 
\begin{figure}
\centering\epsfig{file=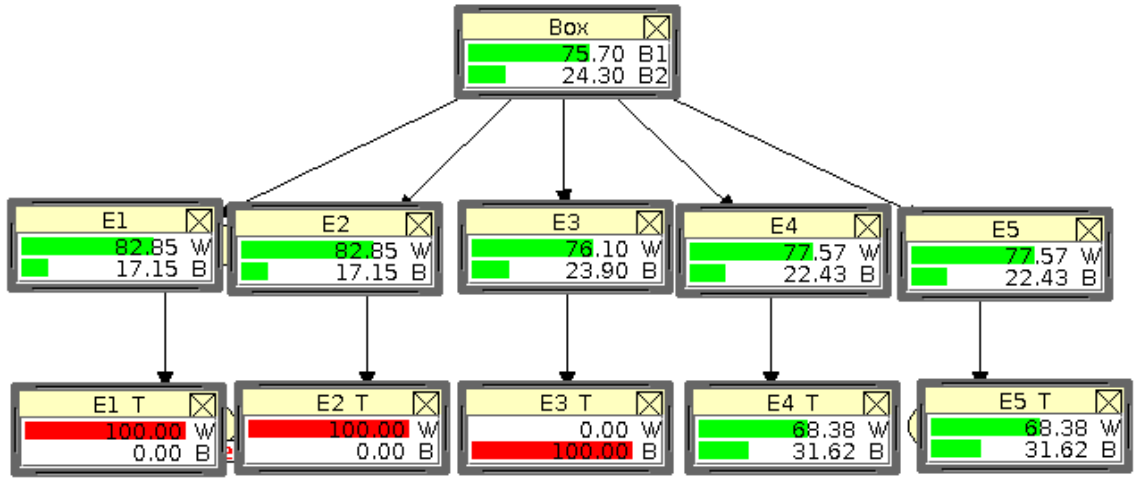,clip=,}
\caption{\small \sf  Network of figure \ref{fig:bn_mon_ET_ww} updated
by a third testimony in favor of black.}
\label{fig:bn_mon_ET_ww_ET_b}
\end{figure}

Let's go on with the experiment and suppose
the third witness says black
(figure \ref{fig:bn_mon_ET_ww_ET_b}). 
This last information reduces the probability of $B_1$, 
but does not falsify this hypothesis, 
as if, instead, we had {\it observed} black.
Obviously, it does also reduce the probability of white balls in the 
following extractions.

The other interesting feature concerns the probability that each 
witness has reported the truth. Our belief that the previous
two witnesses really saw what they said 
is reduced to 83\%. But, nevertheless we are more confident on 
the first two witnesses than on the third one, 
that we trust only at 76\%,
although the lie factor is the same for
all of them. The result is again in agreement with intuition:
if many witnesses state something and fewer say the opposite,
{\it we tend to believe the majority}, if we initially consider
all witnesses equally reliable.
But a Bayesian network
tells us also
how much we have to believe the many more then the fewer. 

Let us do, also in this case the exercise of calculating the 
effective Bayes factor, using however the first formula
in footnote \ref{fb:PEH=0}:
the effective odds $\tilde O_{H}(B_T,I)$ can be written as
\begin{eqnarray}
\tilde O_{B_1}(B_T,I)
&=& \frac{1}{P(W\,|\,B_2)+P(B\,|\,B_2)/\lambda},
\end{eqnarray}
i.e. ${1}/{[1/13 + (12/13)/(1/5)]} = {13}/{61} = 0.213$,
smaller then 1 because they provide an evidence 
against box $B_1$
($\Delta\mbox{JL}=-0.67$). It is also easy to check that the 
resulting probability of 75.7\% of $B_1$ can be obtained 
summing up the three weights of evidence, two in favor of 
$B_1$ and two against it: $\Delta\mbox{JL}_{B_1}(W_T,W_T,B_T,I)= 0.58+0.58-0.67=
0.49$, i.e. $\tilde O_{B_1}(W_T,W_T,B_T,I)=10^{0.49}=3.1$, that
gives a probability of $B_1$ of 3.1/(1+3.1)=76\%.

\begin{figure}[t]
\centering\epsfig{file=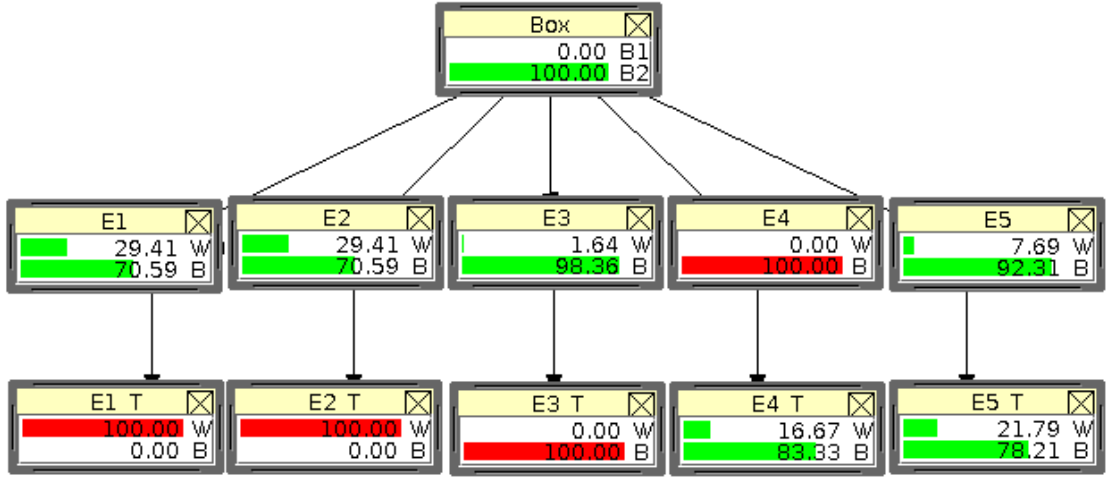,clip=,}
\caption{\small \sf Network of figure \ref{fig:bn_mon_ET_ww_ET_b} updated
by a direct observation of a black ball.}
\label{fig:bn_mon_ET_ww_ET_b_Eb}
\end{figure}

Finally, let us see what happens if we really see a black ball
($E_4$ in figure \ref{fig:bn_mon_ET_ww_ET_b_Eb}).
Only in this case we become certain that the box is of the kind 
$B_2$, and the game is, to say, finished. But, nevertheless,
we still remain in a state on uncertainty with respect to several
things. The first one is the probability of a white ball in future
extractions, that, from now becomes 
1/13, i.e. 7.7\%, and does not change any longer. 
But we also remain uncertain on whether the witnesses 
told us the truth, because what they said is not 
incompatible with the box composition. But, and again in qualitative
agreement with the intuition, we trust much more whom told
black (1.6\% he lied) than the two who told white (70.6\% they lied). 

Another interesting way of analyzing the final network is to consider
the probability of a black ball 
in the five extractions considered. 
The fourth is one, because we have seen it. The fifth is 92.3\% ($12/13$)
because we know the box composition. But in the 
first two extractions the probability
is smaller than it (70.6\%), while in the third is higher (98.4\%). That is 
because in the two different cases 
we had an evidence respectively against and in favor of them.
%

%Let us end with an exercise: try to calculate the probability
%that the witness lied as a function of the available 
%information.
%\mbox{}
%
%\\ \break \mbox{}

\end{document}